\documentclass[10pt,amsfonts,amsmath, fleqn]{article}
\usepackage{multirow}
\usepackage{eurosym}
\usepackage{xcolor}
\usepackage{booktabs}
\usepackage[latin1]{inputenc}
\usepackage{amsmath}
\usepackage{amsfonts}
\usepackage{amssymb}
\usepackage{float}
\usepackage{comment}
\usepackage[ruled,vlined]{algorithm2e}
\usepackage[english]{babel}
\usepackage[margin= 0.7 in, top=0.5 in]{geometry}
\usepackage{amsfonts,graphicx,latexsym}
\usepackage{natbib}
\newcommand{\TUS}{\textit{TUS}}

    \title{A heuristic algorithm for the Air Transport Unit Consolidation Problem}
\author{\it Lorenzo Peirano$ ^{(1)}$ \quad Enrico Angelelli$ ^{(1)}$ \quad Claudia Archetti$ ^{(2)}$\\
	{\small $^{(1)}$ \it University of Brescia, Department of Economics and Management, Brescia, Italy}\\
	{\small lorenzo.peirano@unibs.it,enrico.angelelli@unibs.it}\\
	{\small $^{(2)}$ \it ESSEC Business School in Paris, Department of Information Systems, Decision Sciences and Statistics, Cergy, France}\\
	{\small archetti@essec.edu}}
\date{}

\begin{document}
	
	\maketitle
	
	\begin{abstract}
Consolidation of loose packages into transport units is a fundamental activity offered by logistics service providers. Moving the transport units instead of loose packages is faster (with one movement only, multiple packages are loaded instead of having one load operation for each package), safer (chances of damage and loss is reduced) and cheaper. One of the typical objective of consolidation problems is the minimization of the number of transport units used, e.g. containers. In air transportation, however, transport units have multiple aspects which concur in the calculation of the cost and thus optimization in the number and characteristics of the transport units is required. In this paper, we present the air transport unit consolidation problem where the aim is to determine how to consolidate loose packages in transport units with the goal of minimizing the corresponding cost. The problem is a variant of the three-dimensional bin packing problem where the objective function is formulated according to the way costs are calculated in the air transport business. In addition, side constraints are included to take into account specific requirements. We propose a heuristic algorithm which construct an initial feasible solution and then improves it through a local search algorithm. Computational tests on randomly generated instances show that the algorithm provides high quality solutions in a reasonable computing time. In addition, tests on real data show that it improves solutions found in practice by the company providing the data.
	\end{abstract}


\section{Introduction}\label{sec:desc}
Freight forwarders handle international shipments for their customers and focus their activity on offering the most complete logistic service from origin to destination. This involves not only the management of documents and the physical transfer of goods, but includes other logistic operations (packaging and labelling are just some examples).
One of the activities carried out is the consolidation of loose packages into transport units (TU). For air shipments TUs are usually pallets of various dimensions, and crates. Since most of the times shipper companies are not able to consolidate themselves the goods of which the shipment is composed, freight forwarders spend a high amount of time in performing this operation. The most evident reason behind consolidation in TU is that the movement of a TU is not only faster, more efficient and cheaper than handling a higher number of packages, it also increases the safety of the shipment since it is more difficult to loose or damage packages. Moreover, since the information included in the documentation must be aligned to the characteristics of the shipment, controls by any monitoring organization (customs, FDA, police, etc.) are easier, and, thus, faster. However, choosing the number and quality of TU to be used is not a simple task and it requires the consideration of different aspects. Indeed, the way in which loose packages are consolidated impacts in the shipment cost.
In fact, TUs can have different dimensions and weight capacities, which generates different costs. In the following, we assume the  measures are expressed in centimeters and in the \textit{width $\times$ length $\times$ height} format, with the latter sometimes not reported. The most common TU used is pallets. EUR and EUR1 pallets are the standard European 120x80x15 pallets, EUR2 are 120x100 while EUR6 are 80x60. Height is the same for every pallet and from now on it will not be reported anymore. American standards are slightly different (40x48 inches) while in Japan the standard measure is 110x110.\\
Air transportation services apply a unitary freight on the taxable weight, which is the higher between the real weight of the shipment and the weight equivalent conversion of its volume. Every cubic meter is equal to approximately 167 kilograms of weight.
When the TU is stackable, the volumetric weight is converted on the actual volume of the TU. On the other side, when the TU is not stackable, it effectively occupy the volume generated by the projection of the base dimensions to the maximum height of the vehicle.\\
Each airplane model has different sizes of the cargo hull and thus its limitations and characteristics must be considered. Air transportation is carried out by two main \text{families} of aircraft: Passengers and Freighter. Passenger aircrafts (PAX) transport passengers on the upper deck and cargo on the lower section of the aircraft. Higher priority is given to passenger's luggages and the remaining slots are allocated to cargo shipments. Depending on the model of the aircraft, different kind of Load Unit (LU) can fit inside the lower deck, defining the constraints in loading the cargo into the plane. The most common internal heights for LU for passenger airplane are 130 and 160 cm for medium and long range flights, respectively, while short range flight usually requires smaller airplanes with maximum height on the lower deck of 110 cm. Very small aircrafts do not use LUs or any kind of consolidated shipment and only accept loose packages with very limited weight and dimensions.\\
LUs shape and size defines the maximum dimensions allowed for each TU to be accepted on a given airplane. The TU is limited in every dimension by the physical boundaries of the LU. However, when using standardized TU like pallets, their base is always loadable within any LU and the only dimension which is limited is the height.\\
Freighter (Cargo Aircraft Only, or CAO) aircraft do not transport passengers and thus can make use of the upper section deck to transport cargo. The upper section deck exploits the circular shape of the aircraft to allow the transportation of bigger packages. Most of the times the upper deck do not use LUs and directly load the consolidated TU. One strict condition imposed by air companies is that the material to be loaded on the upper deck must be forkliftable and movable with handling vehicles.\\

\begin{figure}[H]
	\caption{Cross-section of a CAO}
	\includegraphics[scale=0.40]{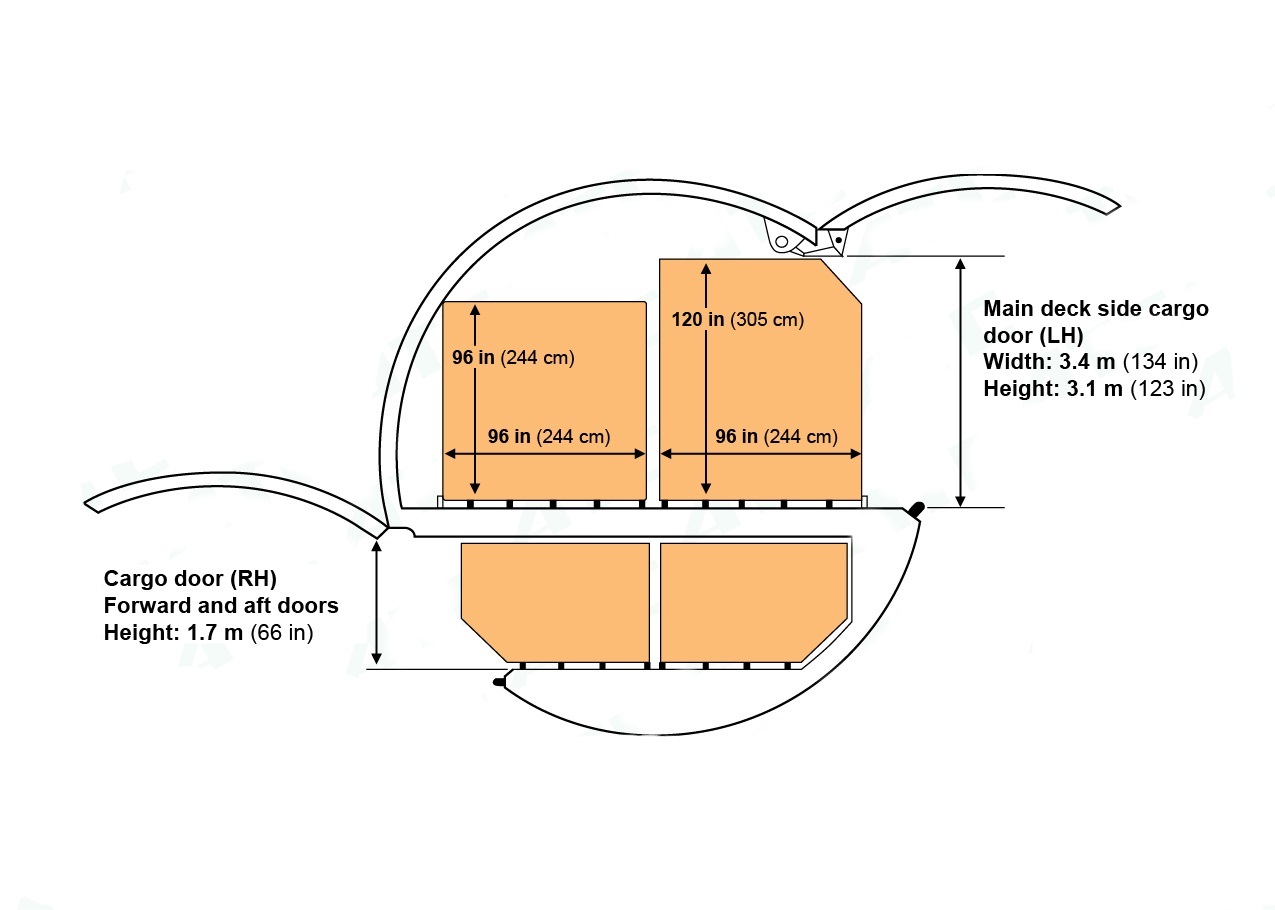}\label{fig:Cross-section}
	\centering
\end{figure}

In Figure \ref{fig:Cross-section} a cross-section of a CAO aircraft is shown. The lower part is the lower deck and it is loaded with LU only. Upper deck is loaded with both LU and loose packages, especially if their dimensions are out of the gauge for the LU. The upper deck height limit depend on the base dimension of the package: with a small base it can be placed near the center where the shape of the aircraft allows for a higher height, otherwise, the wider the base, the lower is the maximum height allowed. For aircrafts with a side door the height of the door is usually considered as maximum height allowed for transport.
CAO aircrafts have a maximum loading height for cargo of 290 cm for medium sized aircraft and 310-320 for large sized aircraft. Even bigger CAO aircraft exists but do not operate with a line service and rely purely on chartering and thus are not considered in this work.

\begin{figure}[H]
	\caption{Boxes, transport units and load units}
	\includegraphics[scale=0.50]{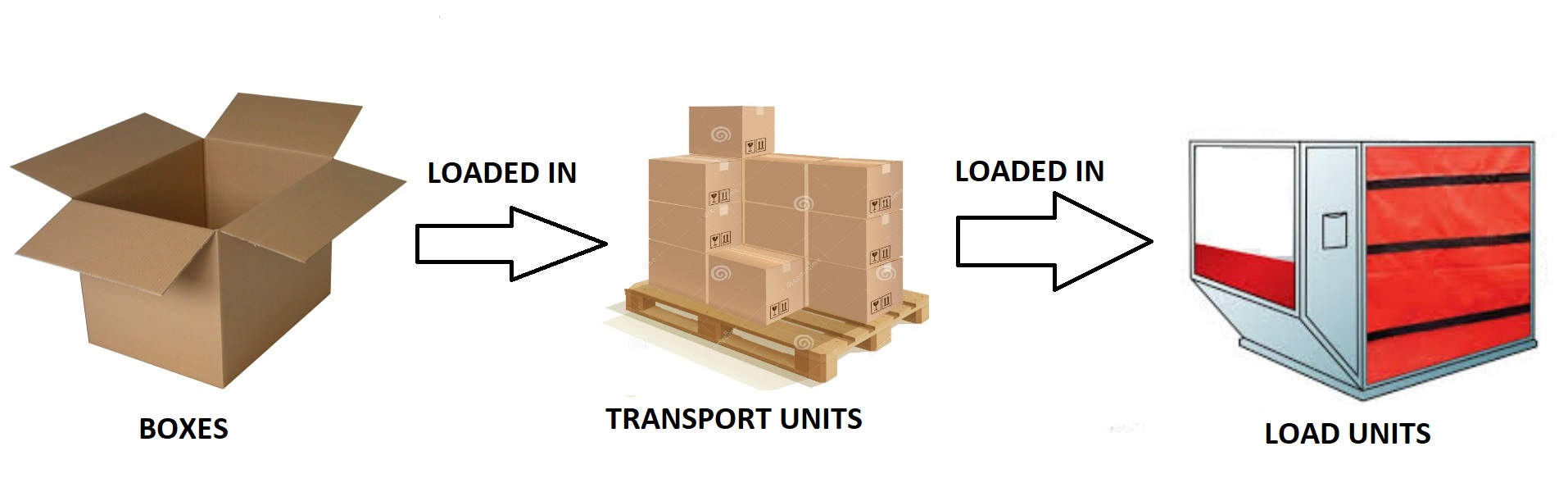}\label{fig:box_tu_lu}
	\centering
\end{figure}
In this paper we study the problem of consolidating loose packages (boxes) in TU to be transported through air transportation services. As shown in Figure \ref{fig:box_tu_lu}, we want to consolidate boxes into pallets, pallets are themself then loaded into LUs by the air company (or by a third party logistics company appointed by the air company). We focus on the first part of the problem, i.e., consolidating loose packages in TUs and we refer to \citet{paquay2018tailored} for a reference on the problem of loading already palletized items or loose packages into LU suitable for air transportation. The problem tackled in \citet{paquay2018tailored} focuses on loading loose packages and palletized packages into LUs. This operation is carried out by air companies (or, on their behalf, by airport handling operators) and happens after the problem we want to tackle in this work. Here, we focus on the first step where a freight forwarder (or a shipper) need to consolidate multiple packages into TU of multiple and varied dimensions in order to minimize the costs that will be paid for the transportation by air means. This is a separate and distinct operation carried out by shippers / freight forwarders on a previous time before delivering the shipment to the airport. While the procedure can be considered as a whole singular operation when shippers/freight forwarder have the availability of LUs, this is not a very common scenario. Most of the time they deliver loose packages or packages consolidated in TUs, and therefore the first operation illustrated in Figure \ref{fig:box_tu_lu} need to be tackled separately.\\
The Air Transport Unit Consolidation Problem (ATUCP) is thus the problem where, given a set of loose packages, the best packing into TUs has to be determined in order to minimize the total cost.
Based on the definition given in \citet{bortfeldt2013constraints}, the ATUCP can be considered as a Multiple Bin-Size Bin Packing Problem (MBSBPP) since the boxes to be packed are extremely heterogeneous, while pallets are weakly heterogeneous. We want to consolidate a set of boxes into the optimal heterogeneous mix of TUs of any given type that seeks the minimization of an objective function considering the value of the TUs (i.e. their taxable weight) and the position of the center of gravity. As explained later, the latter term measures the stability of the TU: we want a TU as stable as possible to avoid damages during transportation. The taxable weight is the unit of measure to which the unitary freight is applied by air companies. Clearly, it cannot be lower than the sum of the weights of the boxes. When weight is low, though, the volumetric weight conversion is usually higher and thus a proper consolidating technique is needed to effectively reduce the cost of the shipment.\\
Boxes can be stackable or non stackable. In the following we treat the case where TU are always non-stackable, i.e., no TU can be loaded on top of another TU. Whenever a TU is non stackable, the volume is calculated applying the maximum height of the airplane deck. Usually this maximum is considered using a standard height of 160 cm for lower deck cargo and 290 cm for upper deck cargo. Moreover, boxes can or cannot be turnable, i.e. two dimensions can be swapped. It is assumed that base dimensions are always swappable, the length can become the width and vice versa. This is not always the case if other means of transportation are considered (for instance, if we consider the consolidation of goods to be shipped through containerization, a package can have extraordinary dimensions that allow it to be forkliftable  on one side only, and thus it can be loaded only with the dimensions derived by the given rotation).\\
One aspect to be taken into account when planning the load layout of a TU is the position of the center of gravity of the load. The center of gravity is the unique point where the weighted relative position of the distributed mass sums up to zero.  This is the point to which a force may be applied to cause a linear acceleration without an angular acceleration. Assuming that the weight is always equally distributed, every box have its center of gravity equal in geometrical center, and the center of gravity of the TU is the weighted average point of the boxes' centers of gravity. A well balanced TU with a center of gravity placed near the central point of the base is easier to transport when lifted with forklifts or cranes, is more stable due to a better balance and, especially when pallets are considered, is more resistant and less subject to damages to the base. Usually it is better to keep heavy items at the bottom of the TU in order to avoid damages to other boxes and improve the balance of the cargo when subject to horizontal movement. Therefore, one of the objective we seek is to place the center of gravity as near as possible to the center of the base and in the lowest position possible. Thus, the objective of the ATUCP combines the minimization of total costs and the optimization of the center of gravity. A numerical example of cost calculation is shown in Section \ref{sec:cost_example} while the calculation of the center of gravity is explained in Section \ref{sec:prob_form}.\\

The contribution of this paper can be summarized as follows:
\begin{enumerate}
	\item We introduced the Air Transport Unit Consolidation Problem (ATUCP) which is a problem finding relevant practical applications in airline freight transportation.
	\item We consider an objective function including solution cost, measured on the basis of the TUs used to consolidate packages, and stability, measured by the center of gravity which depends on the disposition of packages in each TU.
	\item We present a heuristic algorithm which first construct a feasible solution and then improves it through a local search based on different neighbourhoods.
	\item We propose a procedure for generating ATUCP instances which enables a comparison between the solution provided by the heuristic algorithm and an optimal solution determined analytically.
	\item We perform exhaustive tests on randomly generated instances and on real data. The results show that the heuristic algorithm provides high-quality solutions and largely outperforms solutions generated by the operators of the company providing the real data.
\end{enumerate}

The paper in organized as follows. Section \ref{sec:cost_example} presents a numerical example of how solution cost is calculated. The relevant literature is reviewed in Section \ref{sec:literature} while the problem is formally defined in Section \ref{sec:prob_descr}. Section \ref{s:sol_algo} is devoted to the description of the algorithm while computational results are presented in Section \ref{sec:res}. Finally, some conclusions are drawn in Section \ref{sec:conc}.

\subsection{Example of cost calculation}\label{sec:cost_example}

Consider the following simple example: we need to transport 12 boxes with dimensions 60x40x60 and we have availability of pallets with base 120x80. Assume that the weight of boxes is very limited and therefore only the volume is considered. We can choose to consolidate the boxes in one single pallet 120x80x195 (solution 1) or two pallets where one has dimensions 120x80x135 and the other has dimensions 120x80x75 (solution 2).\\
The following table report the volumetric weight for the two possible consolidations.

\begin{table}[H]
    \centering
\caption{Weights} \label{t:pesi1}
%\resizebox*{0.35\textwidth}{!}{
	\begin{tabular}[H]{r|c}
		 & Volume Conversion  \\ \hline
		Solution 1  & 465 kg \\ \hline
		Solution 2  & 513 kg \\

		  \end{tabular}
%	  }
\end{table}

Consolidating in one pallet generates less volume ($120\cdot80\cdot290 = 2.784 $ CBM$ \cdot 167 = 465 $ kg ), but forces to fly with a freighter aircraft since its height is greater than the maximum height of the LU (160 cm). Pallets with lower height on the other side generates more taxable weight ($120 \cdot 80 \cdot 160 \cdot 2 = 3.072 $ CBM $\cdot 167 = 513$ kg). Cargo flights are typically less frequent and more expensive. For instance, if a passenger airline (PAX) company offers a freight of \euro \ 2.00/kg and a cargo airline (CAO) company offers a freight of \euro \ 2.30/kg, the total freight paid is as per table \ref{t:costi1} below.

\begin{table}[H]
	\centering
	\caption{Costs} \label{t:costi1}
	\resizebox*{0.25\textwidth}{!}{
		\begin{tabular}[H]{r|c}
			& Total Cost   \\ \hline
			Solution 1  & \euro \ 1,069.50 \\ \hline
			Solution 2  & \euro  \ 1,026.00 \\

		\end{tabular}
	}
	\end{table}

In this case, creating two pallets is the cheapest option, since it reduces the total cost of the air shipment.
Furthermore, if we consider the availability of pallets with base 120x120 we can consolidate all the boxes on a single pallet 120x120x135 with a volume weight conversion equal to 384 kg.
Considering the freight of PAX flight \euro \ 2.00 the cost is \euro \ 768.00. Consolidating the boxes in a single pallet with base 120x120 is then the best option. In our work it is assumed that there always is a wide choice of flights and aircraft, thus the need for optimization. While during the COVID-19 crisis the number of passenger's flight reduced significantly (even though many of them were used to transport cargo nontheless), the situation returned to pre-pandemic levels, where multiple companies offer their services for the transportation of cargo, using both passenger's and Cargo Aircraft Only options. The effective availability of space on aircrafts impacts the time needed to transport goods and the cost of it, but it is not depending on the consolidation layout and is not directly handled by shippers/forwarders and thus, it is not considered in this study.

\subsection{Related literature}\label{sec:literature}

The ATUCP is a three dimensional bin packing problem (3D-BPP) where the objective is finding the optimal way to allocate loose packages to TUs in order to seek the minimization of TU's cost function and optimize the center of gravity. \citet{martello2000three} give a description of the 3D-BPP and propose an exact branch and bound algorithm for the solution of large instances. A MILP formulation with different valid inequalities is presented by \citet{hifi2010linear}. 3D-BPPs are often very difficult to solve with exact methods. Henceforth, a wide range of heuristics approaches has been introduced in the literature. Tabu search techniques are presented in \citet{lodi2002heuristic} and \citet{crainic2009ts2pack} with good results. \citet{gendreau2006tabu} propose a tabu search algorithm capable of solving a combined capacitated vehicle routing and loading problem. \citet{crainic2009ts2pack} propose a two level tabu-search heuristic. The first level heuristic is meant to deal with the optimality of the bin packing problem, while the second level of the heuristic finds feasible layouts for the packing of the items in the bins. Neighborhoods are explored using a k-chain-moves procedure involving chains of $k$ changes in the solution, which proves to be effective in exploring large neighborhoods in contained computational times. This is clear in the computational results where the method proposed result particularly effective in finding good solution when the time available is limited. Similarly, the three-dimensional single bin-size bin packing problem is studied in \citet{hifi2010heuristics} where the authors present a MILP formulation of the problem. Moreover, they propose valid inequalities usefull to improve the relaxed lower bound of the formulation proposed. Final solutions are found through an heuristic approach. \citet{mack2012heuristic} present a straightforward heuristic capable of solving two-dimensional and three-dimensional bin packing problems. In the same work the authors introduce a set of 1,800 benchmark instances with on average 1,500 items each, which are useful in order to test bin pacing algorithms on very large instances. \citet{de2010particular} present an heuristic for a particular case of the container loading problem where special constraints on the position of the boxes are implemented. The heuristic approach proposed by the author is designed to be easy to implement for warehouse workers during the loading operation and the results of the computational campaign shows that the methodology is capable of generating good solutions, especially when real-world data is used as input. \citet{alvarez2013grasp} present a two phase algorithm useful for solving two and three-dimensional multiple size bin packing problems. The procedure proposed first seek a starting solution using a GRASP heuristic including a constructive phase, a postprocessing phase and improvement moves. The best solution found by the GRASP heuristic are then combined together using a Path Relinking procedure, proposed in three different versions. They compare the computational results with literature benchmark instances. In their work \citet{trivella2016load} present an application of the container loading problem with the additional practical extension of considering the load balance. This latter objective consist in finding the average center of mass in the most suitable position. First, the authors introduce a MILP formulation for the problem. Given the difficulties encountered to solve even small instances, a multi-level local search heuristic approach is presented. \citet{alvarez2015lower} present a study with the aim of improving the lower bounds for integer programming formulations for some relaxations of the original problem. Logical considerations, such as avoiding symmetrical solutions, fixing boxes into bins and setting bound for each type of bins, are implemented in order to improve the bounds for the problem. An exact MILP formulation designed to accomodate packages by layers is designed by \citet{taylor2017three}. The authors show that when the data is pre-processed in order to arrange the order of the boxes to be packed by non-increasing height or volume, calculating very tight lower bounds is an operation of little complexity, which allow the formulation proposed to be solved in reasonable times even for large instances. In general, the objective of the pallet loading problem is the minimization of the number of pallets used to pack the loose packages available, while ensuring their stability. In the ATUCP the object is slightly different and minimizing the number of pallets it's not always the best option. The minimization of the volume (and thus the cost generate for shipping by air means) is part of the objective as well, considering th multiple dimensions of TUs available, it is possible that having a greater number of smaller TUs proves to be a better solution.\\
\citet{elhedhli2017mip} study the mixed-case palletization problem, which is an extension of the three-dimensional bin packing problem which requires the boxes laying upon other boxed to be properly supported. In order to solve the problem the authors propose a MILP-based slicing heuristic which aims at maximizing the fill rate of the lower levels in order to achieve a higher density, allowing upper layers to be placed upon them. The iterative process pack subsets of similar boxes separately, ideally starting from bigger boxes and then filling the gaps remaining by positioning smaller boxes. The methodology was tested against well-known instances in the literature. Further studies on the mixed-case palletization problem are considered in \citet{elhedhli2019three} where the authors propose a novel formulation and a column generation approach suitable to solve it. The sub-problem to optimize is a two-dimensional layer generation problem making use of the concept of \textit{superitems}, where a superitem is a collection of individual items that are compactly stacked together with no gaps between them. The superitems are used to ensure the compactness of the layers. Computational experiments are carried on both benchmark instances derived by the literature and from instances built using industrial data.
  \\
Genetic algorithms proved to be very successful in solving the 3D-BPP: examples can be found in \citet{wu2010three}, \citet{kang2012hybrid}, \citet{gonccalves2013biased}, \citet{li2014genetic}, \citet{dornas2017real} and \citet{ha2017online}. Approximation algorithms are proposed in \citet{miyazawa2009three}.\\
One interesting application is proposed in \citet{wu2017three}. The authors divide their problem in three different layer. The first one requires to pack irregularly shaped item into boxes in order to maximize the utilization rate of the boxes. Afterwards, the boxes are packed in to crates, where the dimensions of the crates could vary significantly and needs to be optimized, leading to a three-dimensional variable-bin packing problem. Finally, the crates need to be loaded into containers. The latter problem is a more classic single container loading problem. The solution is found through the use of a three-stage heuristics. For the second problem, the heuristic at first pack the boxes into a minimum number of crates of larger dimensions, than try to improve the solution by checking if crates of smaller dimensions are a better option. The effectiveness of the heuristic is tested on real-world instances inspired by an industrial case. \citet{alonso2019mathematical} propose mathematical formulations for problems arising in practical applications in logistics. An industry application evolved from a layer-based column generation approach is presented by \citet{gzara2020pallet}, where real-life constraints widely applied by industries are considered and an extensive experimental campaign is provided.
None of the papers mentioned above tackle the ATUCP. Pallet loading problems and container loading problems are similar in their objective of finding stable, compact, balanced loading that maximize the volume occupation of the container (container loading) or minimize the number of containers used (pallet loading problem). Both the problems, in their classical representation, uses one kind of container/pallet of fixed dimensions. In the ATUPC, we consider multiple dimension for the TU that can be used, a highly heterogeneous set of boxes, and the objective we want to optimize takes into account multiple aspects that are not considered together in previous works. The main contribution of this work is the introduction of such a unique problem. Thus, we now provide a formal description of the problem.

\section{Problem Description}\label{sec:prob_descr}

In this section we introduce the notation and the formal definition of the problem.

\subsection{Problem setting and notation}\label{sec:prob_setting}

\subsubsection{Transportation units}
Several types of TU may be available for consolidation, we denote as \TUS\ the set of different types of TU at hand.
Each TU type $tut \in \TUS$ is associated with the following parameters:
\begin{itemize}
	\item $X_{tut}$: the TU width;
	\item $Y_{tut}$: the TU length;
	\item $Z_{tut}$: the TU maximum height allowed;
	\item $Q_{tut}$: the maximum weight that can be carried by the TU;
	\item $V_{tut}=X_{tut}\cdot Y_{tut}\cdot Z_{tut}$: the maximum box volume that can be carried by the TU;
%	\item $CG_{tu}$: Center of gravity of the TU. It corresponds to a point and it is determined by its coordinates within the TU $x_{CG}, y_{CG}, z_{CG}$.
	
\end{itemize}

As procurement of TU items in general is not an issue, we assume to have an unlimited availability for each type of TU $tut \in \TUS$.
In the following we will use $tu \in TU$ to refer to a generic TU item of any type and we will denote as $X_{tu}$, $Y_{tu}$, $Z_{tu}$, $Q_{tu}$ and $V_{tu}$ the respective measures inherited from its type $tut \in \TUS$.

\begin{figure}[H]
\caption{Coordinate system}
\includegraphics[scale=0.40]{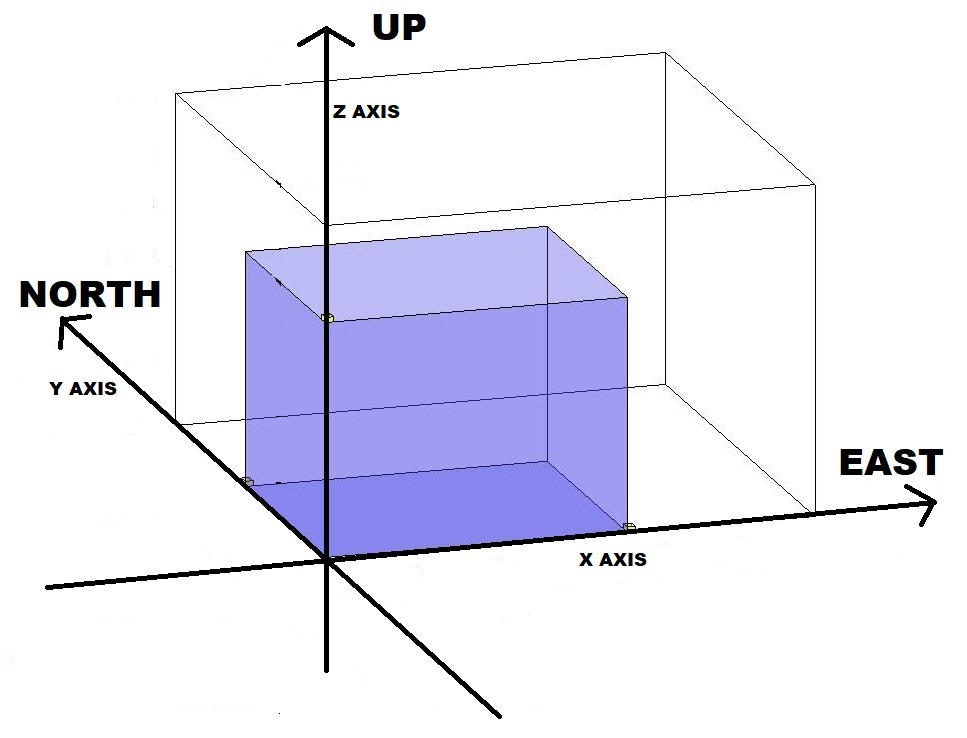}\label{fig:coordinates}
\centering
\end{figure}

Position of consolidated boxes on each TU will be described according to a three-axes Cartesian coordinate system with origin in the south-west corner of the TU base.
In other words, we consider the base as a geographical map, so we move along the $X$ axis from west to east, along the $Y$ axis from south to north, while axis $Z$ starts from the down position and ascends towards the up position.
Thus, the base of a TU $tu$ cover the rectangle $\left[0,X_{tu}\right]\times\left[0,Y_{tu}\right]$ and its center is point $C_{tu}=\frac{1}{2}(X_{tu},Y_{tu},0)$.
The available space for box placement is the rectangular cuboid $\left[0,X_{tu}\right]\times\left[0,Y_{tu}\right]\times\left[0,Z_{tu}\right]$ with volume $V_{tu}=X_{tu}\cdot Y_{tu}\cdot Z_{tu}$.
Figure \ref{fig:coordinates} shows a representation of the coordinate system used.

\subsubsection{Boxes}
\begin{figure}[H]
	\caption{Box orientations}
	\includegraphics[scale=0.35]{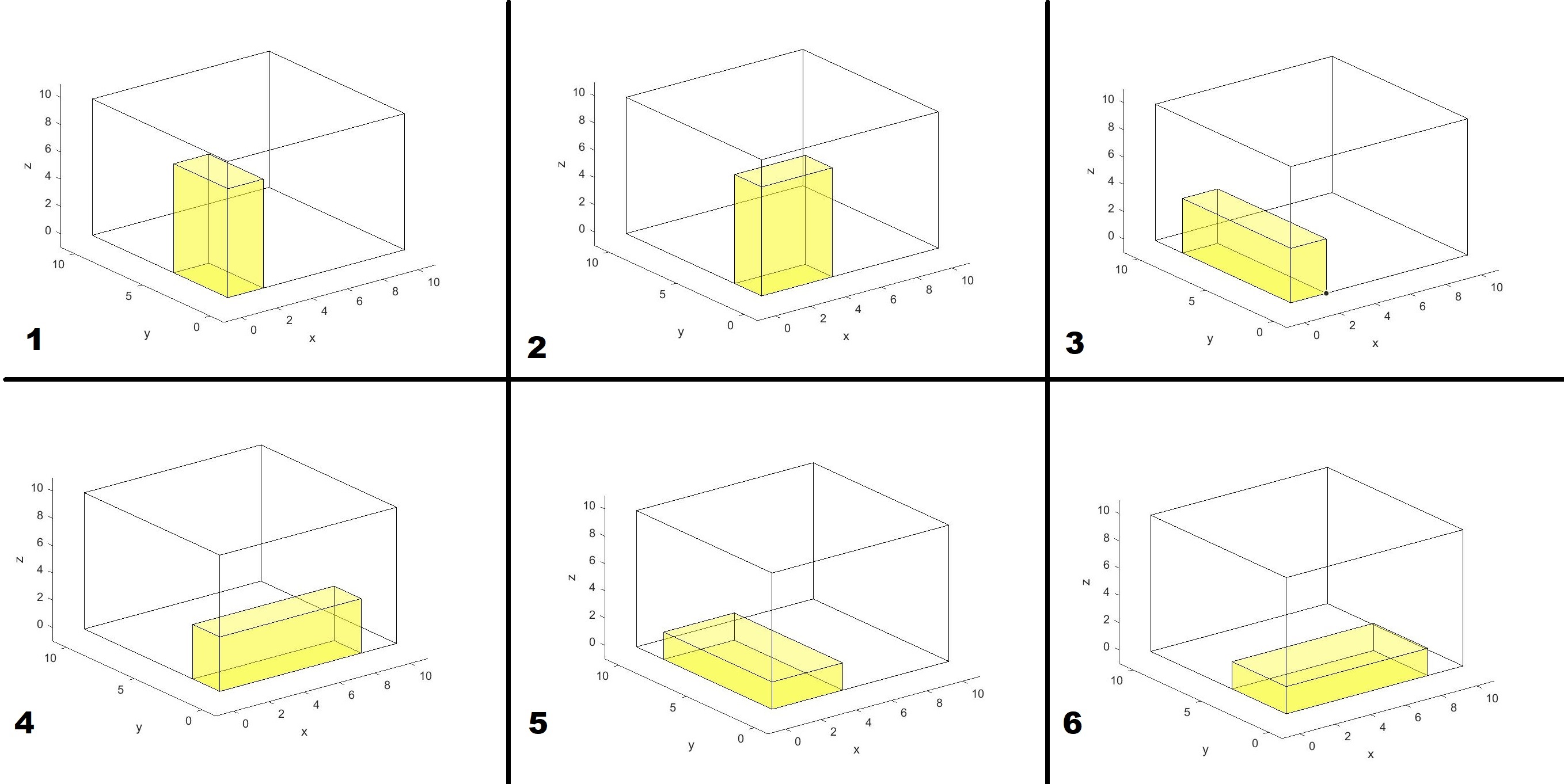}\label{fig:orientations}
	\centering
\end{figure}

We indicate as $B$ the set of boxes to be consolidated.
Each box $b\in B$ has a weight $w_b$ and three dimensions (width, length and height)
that can be potentially oriented in space in six different manners: we have to choose which dimension out of three is parallel to the $Z$ axis, then which of the other twos is parallel to $X$ axis; the remaining dimension is thus forced to be parallel to $Y$ axis.
Flipping dimension orientation (up-down, west-east or north-south) is  irrelevant as this does not change how the box occupies space.
At least two orientations are guaranteed for each box: height parallel to $Z$, and width parallel to either $X$ or $Y$.
Logical attributes $\textit{TXZ}_b$ and $\textit{TYZ}_b$ of box $b$ say whether the box can be oriented with dimension width and length parallel to the $Z$ axis, respectively.
The values of these attributes depend on the characteristic of box content and are fixed for each box.
For sake of simplicity, whatever orientation a box $b$ is given, we denote with $\textit{width}_b$, $\textit{length}_b$, $\textit{height}_b$ its extension along the $X$, $Y$ and $Z$ axes. An example of the six possible orientations is given in figure \ref{fig:orientations}. 
Finally, the logical attribute $st_b$ denotes whether the box is stackable or not. If a box is non-stackable, no boxes can lay upon it.

\subsubsection{Boxes loaded on a TU}
When a box $b$ is assigned to a TU $tu$ we say, with a little abuse of notation, $b \in tu$.

The position of a box $b \in tu$ is described in terms of the coordinate system of $tu$ and we denote as $(x_b,y_b,z_b)$ the position of its west-south-down point in this coordinate system.
As a consequence the box will lay in the 3D interval $[x_b,x_b+\textit{width}_b] \times [y_b,y_b+\textit{length}_b] \times [z_b,z_b+\textit{height}_b]$, and assuming a homogeneous density for the box, its center of gravity will be in point $CG_b \equiv (x_b^{CG}, y_b^{CG}, z_b^{CG})=(x_b, y_b, z_b) + \frac{1}{2}(\textit{width}_b$, $\textit{length}_b$, $\textit{height}_b)$.

The center of gravity $CG_{tu}$ of a loaded TU $tu$  is the weighted position of the centers of gravity of the boxes $b \in tu$. In formula $CG_{tu}\equiv(X_{tu}^{CG},Y_{tu}^{CG},Z_{tu}^{CG}) = \frac{\sum_{b \in tu} CG_b \cdot w_b}{\sum_{b \in tu} w_b} = \frac{\sum_{b \in tu} (x_b^{CG}, y_b^{CG}, z_b^{CG}) \cdot w_b}{\sum_{b \in tu} w_b}$;

\begin{figure}[H]\caption{level curves of $M_{tu}^{xy}$} 		\includegraphics[scale=0.55]{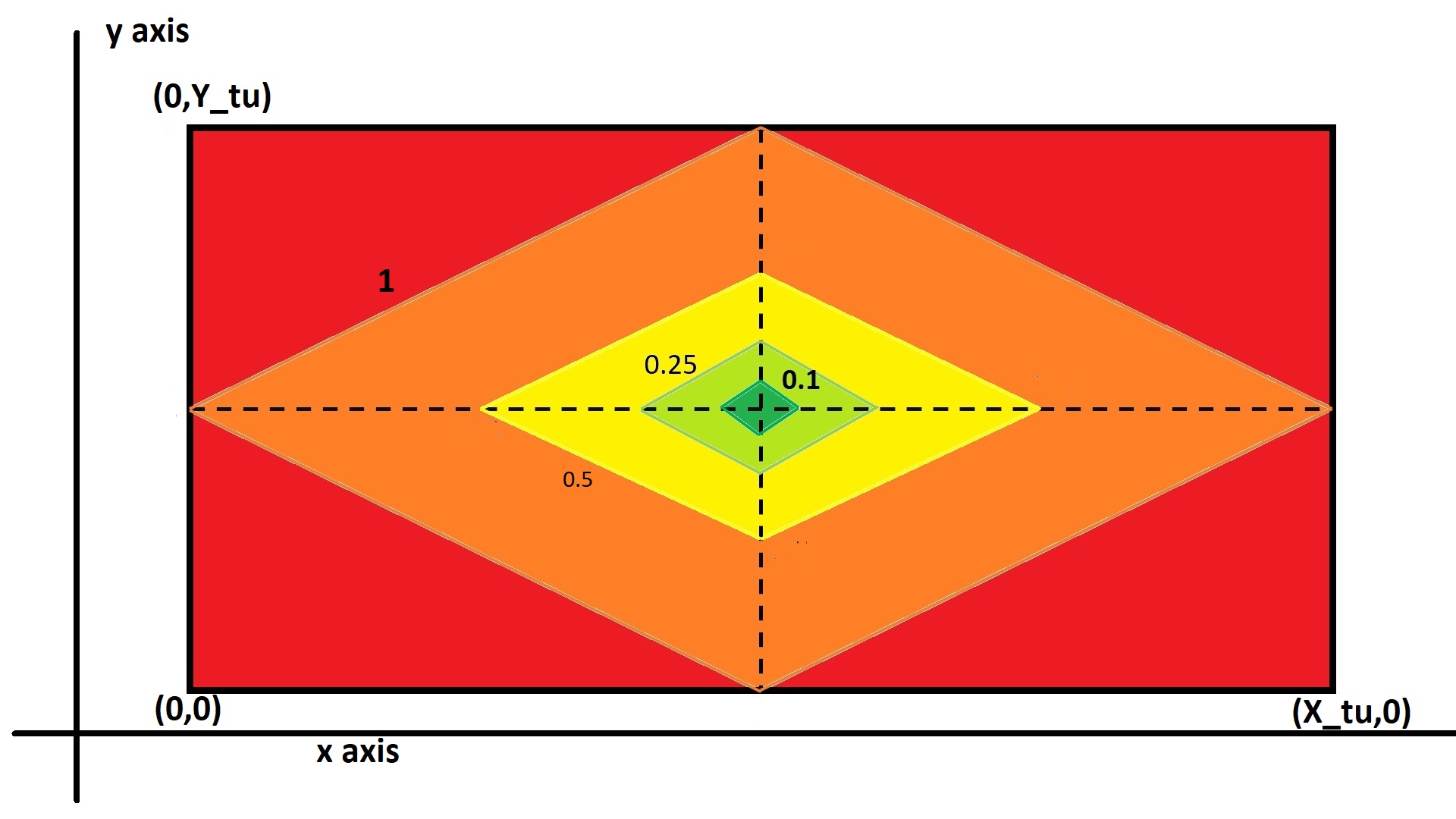}\label{fig:CIxy}		\centering \end{figure}

We also introduce three normalized centering measure of $CG_{tu}$ with respect to the center of the $tu$ base $C_{tu}$, defined as
\begin{eqnarray*}
\nonumber
  M_{tu}^x &=& \left|\frac{2X_{tu}^{CG}-X_{tu}}{X_{tu}}\right| \\
  M_{tu}^y &=& \left|\frac{2Y_{tu}^{CG}-Y_{tu}}{Y_{tu}}\right| \\
  M_{tu}^z &=& \left|\frac{Z_{tu}^{CG}}{Z_{tu}}\right|
\end{eqnarray*}
In particular, measure $M_{tu}^{xy}=M_{tu}^x+M_{tu}^y$ determines how much centered is the projection of $CG_{tu}$ on the base of the $tu$. Ideally, the perfect layout has $M_{tu}^{xy}=0$.
Figure \ref{fig:CIxy} shows level curves of $M_{tu}^{xy}$.

\subsection{Problem formulation}\label{sec:prob_form}
The objective of the ATUCP is to find the feasible consolidation layout of all boxes in $B$ which minimize the value of objective function (\ref{for:OF}), which we call \textit{FITNESS}:

\begin{equation} \label{for:OF}
 FITNESS = \sum_{tu \in TUL} \left ( V_{tu} + \alpha \cdot (M^{xy}_{tu} \cdot  M^{z}_{tu} + \Theta) \right ) + \beta \cdot |TUL|
\end{equation}

where $TUL$ is the list of TUs actually used in the consolidation:

\begin{itemize}
\item $V_{tu}$  is the maximum volume allowed to the TU $tu$. As explained in Section \ref{sec:cost_example}, we use the volume as an approximation of the cost;
\item $M^{xy}_{tu}$ is the already mentioned measure of the distance between the projection of the center of gravity $CG_{tu}$ of the load on the base of $tu$ and its center $C_{tu}$. Ideally, the perfect layout has $M_{tu}^{xy} = 0$ with a balanced distribution of load around the center of the TU.
\item $M^z_{tu}$ is used to evaluate the vertical position of the center of gravity of the TU. Ideally, heavier boxes should be positioned on the bottom, which leads the load center of gravity to lower positions.
\item $\Theta$ and $\alpha$ are tuning parameters aimed at balancing the impact of $CG_{tu}$ and $V_{tu}$ on the objective function. In particular, $\Theta$ is used as a scaling factor to balance the term related to the center of gravity evaluation with the other two terms in the objective function.
\item The term $\lvert TUL \rvert$ is the number of TU used in the consolidation while parameter $\beta$ represents the fixed cost for each TU expressed in volume. Indeed, each used TU generates handling time cost as x-rays check, transshipment, and customs and security inspections.
\end{itemize}

Thus, the objective of the problem is to minimize the number of TUs used and the total volume generated while keeping the center of gravity of the TUs in a position that is as centered and low as possible. While the major objective is the minimization of the total volume, the other objective ensure the stability of the TUs created and the minimization of the costs derived by a higher number of TUs created. The parameters are tuned in order to give priority to the minimization of the costs (generated by volume and number of TUs), while optimizing the center of gravity of the loading layout.
The feasible solution space is defined by the following set of constraints:

\begin{itemize}
\item Box orientation: according to attributes $TXZ_b$ and $TYZ_b$, orientations with width and/or length dimensions parallel to $Z$ axis may not be allowed to some boxes. However when width or length dimensions are allowed to be parallel to Z axis, the height dimension can be parallel to X or Y axis;
\item Weight limit: the sum of box weights assigned to each $tu \in TUL$ must be lower than the TU weight capacity $Q_{tu}$. That is, $\sum_{b \in tu} w_b \leq Q_{tu}$
\item TU boundaries: all boxes assigned to each $tu \in TUL$ must respect the TU boundaries; that is, for all $b\in tu$ the following conditions must hold: 
     $$\left\{\begin{array}{l}
     0\leq x_b \wedge x_b + \textit{width}_b \leq X_{tu}\\
     0\leq y_b \wedge y_b + \textit{length}_b \leq Y_{tu}\\
     0\leq z_b \wedge z_b + \textit{height}_b \leq Z_{tu}
     \end{array}\right.$$
In words, along each axis the edge of the box must be fully contained within the TU boundary. Namely, regarding X axis, the west coordinate must be positive ($0\leq x_b$) and the east coordinate must be not larger than $X_{tu}$ ($x_b+\textit{width}_b \leq X_{tu}$). Analogously. With respect to Y a Z axes.    
\item Overlapping: boxes must not overlap each other, i.e. one box cannot penetrate within another box; overlapping check is given by the statement: two boxes $a,b \in tu$ overlap each other if and only if the following system is satisfied
    \begin{equation} \label{for:CanFit}
    \left\{\begin{aligned}
     \max(x_a, x_b) < \min(x_a+\textit{width}_a,x_b+\textit{width}_b)\\
     \max(y_a, y_b) < \min(y_a+\textit{length}_a,y_b+\textit{length}_b)\\
     \max(z_a, z_b) < \min(z_a+\textit{height}_a,z_b+\textit{heigth}_b)
    \end{aligned}\right.
    \end{equation}
\item Box stackability: The space above non-stackable boxes must be left empty, this means that no boxes can occupy the space of the projection along the Z axis of the base of the non-stackable box.
\end{itemize}

Note that none of these constraints considers vertical stability. The usual approach to ensure vertical stability is to fill the eventual empty spaces with pluriball/alluminium sheets. The proposed packing procedure focuses on trying to limit as much as possible the creation of gaps and empty spaces and henceforth no hard constraints are implemented to ensure vertical stability.

To the best of our knowledge, there is no version of MBSBPP considering this set of constraints in the literature.

\section{An iterated local search algorithm for the ATUCP}\label{s:sol_algo}

We approach the enormous complexity of the problem with a  heuristic algorithm based on the Iterated Local search (ILS) scheme which we call \textit{ILS-ATUCP algorithm} and whose general logic is depicted in the next lines. 
More details follow in the next sections.

\medskip\noindent\textit{The construction algorithm}

Each of the following phases strongly relay on a constructive algorithm called 3DBP that receives in input a type of TU $tut \in \TUS$ and a set of unassigned boxes. 
The algorithm picks one new TU of the given type and feasibly fills it with unassigned boxes; this step is iterated until the set of unassigned boxes is empty.
Details on algorithm 3DBP are provided in Section \ref{sec:3dbp-alg}.

\medskip\noindent\textit{Initialization phase}

In the first step set \TUS\ is sorted according to increasing values of volume $V_{tut}$  and the
constructive algorithm 3DBP is called on the first TU type in the list and the whole set of unassigned boxes.

\medskip\noindent\textit{LS1 search}

Once the initial solution is constructed, a first order local search (LS1), whose moves are described in Section \ref{sec:LS1}, is applied to possibly improve the packing
by moving boxes between TUs or by destroying some TUs and assigning again the boxes by means of algorithm 3DBP in an attempt to reduce the number of utilized TUs and/or better assign part of boxes to TUs.
The type of TU used by algorithm 3DBP in this stage is the last used to build a solution. Namely, when LS1 is executed immediately after the initialization phase, the smallest volume TU type will be used, when LS1 is executed after the LS2 search, the last TU type considered in LS2 will be used in LS1.
The search proceeds until no improvement is found.\\
The aspect of considering different air services with different air freights is not considered, since the selection of the best air service to be used is a complex problem which consider a wide spectrum of technical aspects. However, given a set of TUs that can be loaded in a given air service, the algorithm, through the use of LS2, will search for the best layout possible allowing to pack the loose boxes in a mixed set of different types of TU. If multiple air services are available, the procedure can be repeated for each service by adjusting the available types of TU between those that can be accepted by the air service considered. In this way the best layout for each service will be given by the algorithm.

\medskip\noindent\textit{LS2 search}

When LS1 fails to improve the incumbent solution, a second order local search (LS2) is started.
A set of TUs, whose loads satisfy certain requirements, are removed from the incumbent solution and constructive algorithm 3DBP is repeatedly called to rebuild a solution, every time with a different type of TU.
Every time an improving solution is found through LS2, the algorithm restarts with LS1.
LS2 is applied until all possible TU types are checked in the same round, i.e. when all the possible TU types are tried with no improvements, the algorithm stops.

The reason why we emphasize the distinction between LS1 and LS2 is that while LS1 is devoted to little adjustments of the solution (exchange of boxes between TUs or replacing a number of TUs with other TUs of the same type), LS2 search is aimed at modifying the solution by diversifying TU types. As we show in the computational analysis, the LS2 search is of fundamental importance in solution improvements proving the role of diversification of TUs type in the load consolidation.

\subsection{Algorithm 3DBP}\label{sec:3dbp-alg}

The construction algorithm receives in input a type of TU $tut\in\TUS$ and a list of unassigned boxes $\bar{B}$. The output produced by the algorithm is a feasible assignment of boxes in $\bar{B}$ to an appropriate number of TU of the given type.

\begin{figure}[H]
	\caption{Spatial system}
	\includegraphics[scale=0.45]{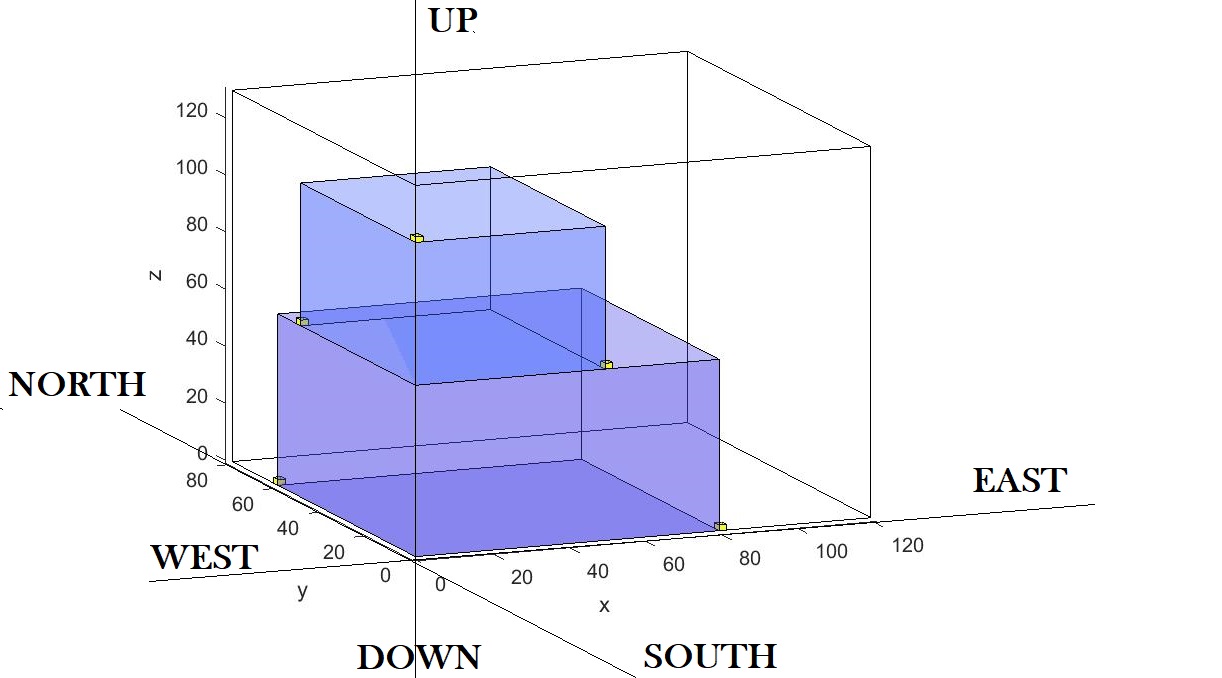}\label{fig:system}
	\centering
\end{figure}

We implemented an adaptation of the Extreme Point (EP) methodology introduced by \citet{crainic2008extreme} for the 3DBP. 
EPs are the positions where we want to place the south-west-down corner of a box. 
For each TU $tu \in \textit{TUL}$, the list $\textit{EPL}_{tu}$ of available EPs is initialized with
the south-west corner of its base
where the origin of the coordinate system associated to the TU is located. Figure \ref{fig:system} show a representation of the spatial system used by the algorithm.
The list is updated each time a new box is added to the TU. 

Every EP $e$ is defined by its coordinates ($x_e$, $y_e$, $z_e$) and three characteristics that also depend on the current load of the TU:

\begin{itemize}
\item $\bar{X}_e$: maximum dimension allowed along the $X$ axis for a box $b$ to be positioned in $e$ (i.e. $width_b \leq \bar{X}_e$). Any box with a larger width would either pass the TU boundary or compenetrate some other box.
\item $\bar{Y}_e$: analogously to $\bar{X}_e$ defines the maximum dimension allowed along the $Y$ axis for a box to be positioned in $e$ (i.e. $length_b \leq \bar{Y}_e$).
\item $\bar{Z}_e$: analogously to $\bar{X}_e$ defines the maximum dimension allowed along the $Z$ axis for a box to be positioned in $e$ (i.e. $height_b \leq \bar{Z}_e$).
\end{itemize}

Given the set $\bar{B}$ of boxes to be packed and one TU type $tut$, the algorithm first sorts the boxes using algorithm $SORT$, described in Section \ref{sec:sort}, to identify the order in which boxes are assigned to TUs.
Once the list of boxes in  $\bar{B}$ is properly sorted, the algorithm starts to pack them one at a time.

Each box is loaded in its best position on one of the TUs already in $TUL$ if possible. 
Otherwise a new TU $tu$ of type $tut$ is added to $TUL$, its EP list $\textit{EPL}_{tu}$ is initialized with $e=(0,0,0)$, and the current box is placed at that EP.

Once the current box is placed, $\textit{EPL}_{tu}$ is updated by using the algorithm EP-Generation described in Section \ref{sec:EP_gen_algo}, the box is removed from $\bar{B}$ and the next box is considered until $\bar{B}$ is empty.

In order to decide the position of a box, for each TU $tu \in TUL$, for every EP $e \in \textit{EPL}_{tu}$ and for every allowed box orientation, algorithm \textit{CanFit} (described in Section \ref{sec:canfit_algo}) checks whether the box can be feasibly placed in $e$ on $tu$.
In this case, the assignment is priced by the formula:
\begin{equation}\label{eq:algo_obj}
\begin{aligned}
\textit{COST}(b,e) = & N \cdot z_e +  x_e +  y_e\\
& + M \cdot (z_e + \textit{height}_b)\\
& - N\cdot\theta\cdot\left[(\bar{X}_e - \textit{width}_b)+ (\bar{Y}_e - \textit{length}_b)\right]\\
& + \lambda \cdot \left[(\bar{X}_e \mod \textit{width}_b) + (\bar{Y}_e \mod \textit{length}_b)\right] - \textit{NBOX}
\end{aligned}
\end{equation}
where $N$ and $M$ are fixed big constants and:
\begin{itemize}
\item $N \cdot z_e +  x_e +  y_e$ is the evaluation given to the position of the EP $e$. Lower first and western-southern EPs are preferred, as this leads to more compact layouts where boxes are more likely to be tight, without empty spaces between them. Even though western-southern points do not explicitly seek for load centering, this guarantee a better exploitation of space in horizontal dimensions.
\item $M \cdot (z_e + \textit{height}_b)$ is the evaluation given to the position of the west-south-up corner of the box. It represents the height of the box within the TU and favors, accordingly with the objective function, orientations that minimizes the height, to keep the center of gravity as low as possible.
\item $N\cdot\theta\cdot\left[(\bar{X}_e - \textit{width}_b) + (\bar{Y}_e - \textit{length}_b)\right]$ represents how effectively the box occupy the remaining space left for the EP along the $X$ and $Y$ axes. Accordingly with the objective function, a better utilization of the empty space leads to more compact layouts and a lower number of TUs.
\item  $\lambda \cdot \left[(\bar{X}_e \mod \textit{width}_b) + (\bar{Y}_e \mod \textit{length}_b)\right]$ is used to place the box in a position that better partition the remaining empty space. By dividing the dimension of the TU by the dimension of the box, if the remainder is zero then boxes of similar dimensions can fit perfectly along that direction. The idea is that when multiple boxes of the same dimensions are available, the algorithm tries to put them oriented in the best possible way to cover the base of the TU. Accordingly with the objective function, a better utilization of the space leads to more compact layouts and a lower number of TUs.
\item $NBOX$ is the number of boxes already placed in the TU. Ideally, each TU should be filled as much as possible. When choosing in which TU a box should be placed, full TUs are preferable to empty TUs, to better utilize the volume of the TUs already created. Accordingly with the objective function, the minimization of this quantity help improve the fitness of the assignment.
\end{itemize}

The algorithm place the current box in a TU and position and orientation which is feasible and with minimum value of function $COST$.

The idea of the pricing formula is to prefer the west-south-down available EP, while taking into account both the height of the box (in order to favor orientations where the height of the package is lower) and the remaining space available after a box has been inserted, to better occupy the empty space of the TU.
Moreover, the pricing formula tries to find a perfect partitioning of the empty space left to obtain tighter layouts, and prefers to place boxes in TUs with the highest number of boxes already placed. 

Function \eqref{eq:algo_obj} components act together in order to find layouts that suits the objectives described in \eqref{for:OF}. 
To summarize, the desired layout minimizes the height of the TUs in $TUL$, thus minimizing the value of $CG_{tu}$. Having compact layouts, where the space utilization is optimized, helps to reduce the total volume $V_{tu}$ and the number of TU used in the final solution. Moreover, compact layouts with no empty spaces within boxes tend naturally to produce centers of gravity toward the center of the TU. Overall, every aspect of the cost function \eqref{eq:algo_obj} is designed to synergize with the general minimization  of objective function \eqref{for:OF}.

When the ideal position of the box is identified, the algorithm inserts the box in the selected TU at the selected EP and:
\begin{itemize}
\item updates the list of available EPs using algorithm EP-Generation (EP occupied by the box is removed and new potential EP points are added to $\textit{EPL}_{tu}$).
\item updates $\bar{X}_e$, $\bar{Y}_e$ and $\bar{Z}_e$ of every EP of the TU.
\end{itemize}

\subsubsection{Algorithm SORT} \label{sec:sort}

Algorithm SORT is focused on creating an ordered list of packages on the basis of weight, dimension and height, as described in the following.
The algorithm is called by algorithm 3DBP an receives in input the type of TU $tut$ algorithm 3DPB is working with, the set of unassigned boxes $\bar{B}$ and two parameters $n$ and $m$  whose role is explained below.

Ideally, heavy boxes needs to be placed at the base of the TU, to optimize the center of gravity objective.
At the same time, big boxes, especially those with a large base, should be allocated first when large space slots are still available while small boxes will be easily placed in a some of many small space slots formed during TU loading operations.
Thus, algorithm SORT starts by rotating each box $b \in \bar{B}$, compatibly with its rotation attributes $\textit{TXZ}_b$ and $\textit{TYZ}_b$ in such a way to obtain the minimum possible height.
Then, the algorithm divides the set $\bar{B}$ in $n$ first order clusters based on box weight.
First order cluster $C(i)$ contains boxes with weight in the range $\frac{Q_{tut}}{n} \cdot\left( (i-1), i\right]\ \forall i=1,\dots,n$.
Afterwards, the algorithm divides each first order cluster in $m$ second order clusters based on box base dimension.
Second order cluster $C(i,j)$ contains boxes from cluster $C(i)$ with base area $(x_b \cdot y_b )$ in the range $\frac{X_{tut} \cdot Y_{tut}}{m} \cdot\left( (j-1), j\right]\ \forall j=1,\dots,m$.
Lastly, each of the $n \cdot m$ clusters is internally sorted with respect to decreasing height of the boxes.

Finally, clusters $C(i,j)$ are taken in decreasing order ($i=n\ldots,1$ and $j=m\ldots,1$ for each $i$) to compose the final sorted list which is sorted by decreasing weight first, decreasing base second, and decreasing height last.

\subsubsection{Algorithm CanFit} \label{sec:canfit_algo}
\begin{figure}[H]
\caption{Boxes placed in the pallet and maximum dimensions for EP (0,0,30)}
\includegraphics[scale=0.60]{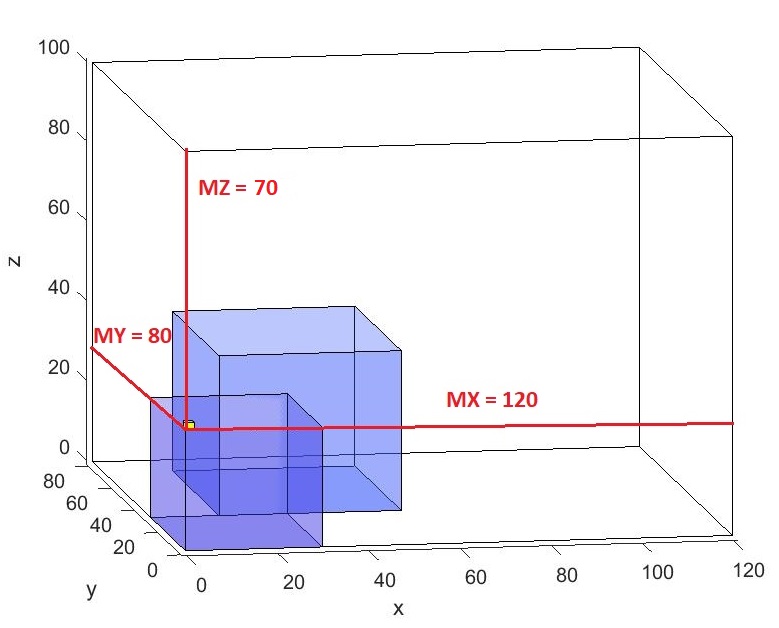}\label{fig:proj}
\centering
\end{figure}

Algorithm CanFit receives in input a TU $tu$ with its current load, a candidate EP $e$ and an oriented box $b$. The algorithm must respond with a yes or no to the inquiry whether box $b$ can be placed on $tu$ in the given EP $e$ with the given orientation.
First, algorithm check whether $b$ dimensions are lower than the maximum dimensions allowed for a package in the EP chosen, that is $width_b \leq \bar{X}_e$, $length_b \leq \bar{Y}_e$, and $height_b \leq \bar{Z}_e$. If just one of the inequalities fails, algorithm stops with a no answer.
Otherwise, we do not have an answer yet and we need to check inequalities \eqref{for:CanFit} where the role of box $a$ is played by each box already loaded on $tu$. If at least one box $a$ exists such that equation \eqref{for:CanFit} is satisfied the algorithm stops with answer no, otherwise return answer yes.
To better illustrate the situation we provide the following example. In a pallet $120\times 80\times 100$ we have, among others, one box $b_1$ with dimension $30\times 30\times 30$ in position (0,0,0) and another box $b_2$  with dimension $40\times 40\times 40$ in position (15,30,0). If we consider EP (0,0,30), the one generated by the south west up corner of the first box,  the maximum dimensions for boxes considered for that EP are $120$ along the $X$ axis, $80$ along the $Y$ axis and $70$ along the $Z$ axis, as showed in Figure \ref{fig:proj}.

\begin{figure}[H]
\caption{Insertion of box $b_3$}
\includegraphics[scale=0.30]{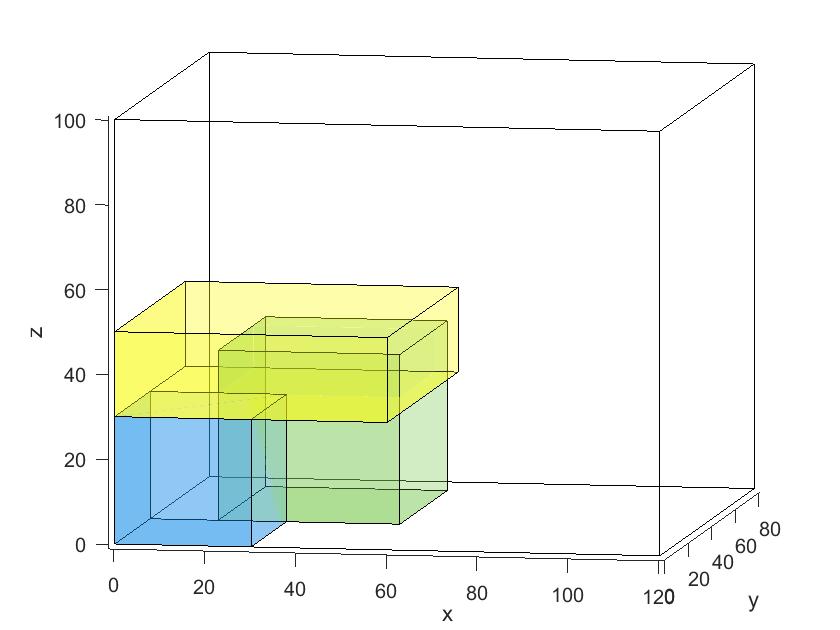}\label{fig:scatoverlappa}
\centering
\end{figure}

Next, we want to insert a box $b_3$ with dimension $60\times 60\times 20$. By applying the first fit check mentioned above, we should be able to insert it in (0,0,30) since the dimensions of the box are lower than the maximum dimensions of the EP. Figure \ref{fig:scatoverlappa} shows the result of this operation.

Clearly, boxes $b_2$ and $b_3$ overlaps. Henceforth, we need a further check based on the position and dimensions of the boxes already placed in the TU in order to avoid overlaps.

\subsubsection{Algorithm EP-Generation} \label{sec:EP_gen_algo}

Each time a box $b$ is allocated to a TU $tu$ the EP $e$ where the box is placed is removed from $\textit{EPL}_{tu}$.
Then up to five new EPs are added to $\textit{EPL}_{tu}$. The new EPs are generated according to the following rules:

\begin{enumerate}
\item The east-south-down corner of box $b$ with coordinates ($x_e + \textit{width}_b$, $y_e$, $z_e$) is projected down, until it touches the base of the TU in the $XY$ plane or the up face of a box. The resulting point is then projected south, until it touches the plane $XZ$ or the north face of a box. The point resulting from the projection is added to $\textit{EPL}_{tu}$.
    More precisely, the projection down either touches the plane $XY$ at point ($x_e + \textit{width}_b$, $y_e$, 0) or it stops on the face of the first package $b_1$ it encounter in point ($x_e + \textit{width}_b$, $y_e$, $z_{e_1}+z_{b_1}$) where $e_1$ is the EP where box $b_1$ is placed. 
    The algorithm proceeds to project the newly found point to the $XZ$ plane until it touches the plane in either point ($x_e + x_b$, $0$, $0$) or point ($x_e + \textit{width}_b$, $0$, $z_{e_1}+z_{b_1}$), or until it touches the face of the first package $b_2$ along the projection in either point ($x_e + \textit{width}_b$, $y_{e_2} + \textit{length}_{b_2}$, 0) or point ($x_e + \textit{width}_b$, $y_{e_2} + \textit{length}_{b_2}$, $z_{e_1}+z_{b_1}$), depending on the previous projection. 
    The point found, either on the face of the TU or on a face of another box, is added to $\textit{EPL}_{tu}$.
\item The west-north-down corner of box $b$ with coordinates ($x_e$, $y_e+\textit{length}_b$, $z_e$) is projected first on the $XY$ plane and then on the $YZ$ plane like in case 1.
\item If box is stackable (attribute $st_b$ is true),  west-south-up corner point with coordinates ($x_e$, $y_e$, $z_e + \textit{height}_b$) is added to $\textit{EPL}_{tu}$.
\item If box is stackable (attribute $st_b$ is true), the west-south-up corner point with coordinates ($x_e$, $y_e$, $z_e + \textit{height}_b$) is projected to the $XZ$ plane until it either touch the plane in point ($x_e$, $0$, $z_e + \textit{height}_b$) or it stops the projection on the face of the first package encountered $b_1$ in point ($x_e$, $y_{e_1} + \textit{length}_{b_1}$, $z_e + \textit{height}_b$). The point resulting from the projection is added to $\textit{EPL}_{tu}$.
\item Similarly to 4), west-south-up corner point with coordinates ($x_e$, $y_e$, $z_e + \textit{height}_b$) is projected to the $YZ$ plane only if the box is stackable. In this case, the projected point is added to $\textit{EPL}_{tu}$.
\end{enumerate}

As a special case we observe that the first box to be assigned to a TU is always placed in $e=(0,0,0)$ and projections (3), (4) and (5) produce the same point $(0,0,\textit{heigth}_b)$, resulting in the generation of 3 points instead of 5.

\begin{figure}[H]
\caption{EP generation example}
\includegraphics[scale=0.75]{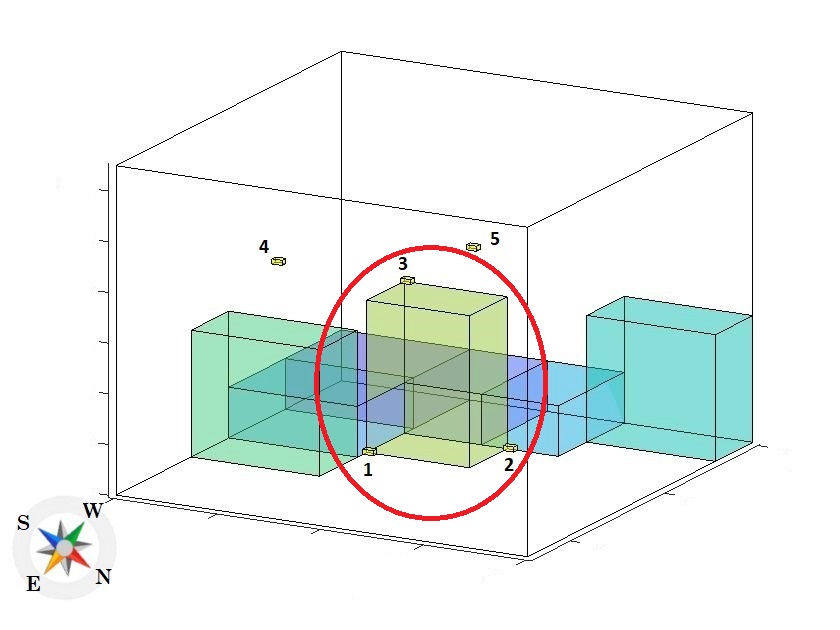}\label{fig:Upd_EP}
\centering
\end{figure}

In Figure \ref{fig:Upd_EP} we show a graphical example of the depicted procedure, where a newly inserted box (circled in red) generates  five EPs. For the sake of clarity, only the newly generated EPs are shown in the picture.

\subsection{First order local search (LS1)}\label{sec:LS1}
Local search LS1 implements three moves $N1$, $N2$ and $N3$ applied sequentially.
If a move succeeds in finding an improving solution, the search is started over from $N1$ with the new incumbent solution. If no move is successful, the search proceeds to LS2.
Implementation of each move consists in applying a determined sequence of steps (micro-moves) for a fixed number of times. Specifically, each step is repeated three times in the implementation whose results are presented in Section \ref{sec:res}.

\begin{itemize}
\item[$N1$:] At each step, one box on top of one TU (origin) is taken and moved on top of another TU (destination). The box to be moved is randomly selected among the three with highest position in the origin TU and placed in the destination TU so that $COST$ function is minimized. Different strategies are used to select which TU is used as origin and which TU is used as destination, i.e. different strategies are employed to decide from which TU move boxes and which boxes move them to. The strategies are employed in the following order:
    \begin{itemize}
    \item the origin TU is the one with the greatest weight; the destination TU is the one with smallest weight;
    \item the origin TU is the one with the greatest weight; the destination TU is randomly chosen;
    \item the origin TU is the one with the greatest height; the destination TU is the one with smallest height;
    \item the origin TU is the one with the greatest height;  the destination TU is randomly chosen;
    \item the origin and destination TUs are randomly chosen;
    \end{itemize}
\item[$N2$:] The iterated sequence of steps is composed as follows:
    \begin{itemize}
    \item two TUs are randomly selected;
    \item one box for each TU is randomly chosen among those that are laying on the top;
    \item if the two boxes can be feasibly swapped their new position is determined so that $COST$ function is minimized for each of the receiving TUs.
    \end{itemize}
\item[$N3$:]The iterated sequence of steps is composed as follows:
    \begin{itemize}
    \item draw a random integer number $n$ between $2$ and $|TUL|$;
    \item destroyed $n$ TUs randomly chosen
    \item call algorithm 3DBP to recover a feasible solution; input for algorithm 3DBP is the set of boxes released by TUs destruction and the last type of TU either used in the initialization step or in the last call from LS2 as explained in Section \ref{sec:LS2}.
    \end{itemize}
\end{itemize}

\subsection{Second order local search (LS2)}\label{sec:LS2}
While Local search LS1 is focused on adjusting load of TUs already packed using multiple moves, LS2 is based on a single type of move. 
The idea is to variate the type of used TU by destroying a set of TUs in the incumbent solution and rebuild a feasible solution by repeated calls to algorithm 3DBP with different types of TU.
Destroyed TUS are those which satisfy at least one of the following criteria:
\begin{itemize}
\item TUs with a filling rate lower than $\Omega$. The idea is that TUs with a low filling rate contain few boxes which might be easily inserted in existing TUs or in TUs of a different type.
\item TUs where the lateral space left is greater than $\gamma$.  If the minimum space left between the boxes and the boundaries of the TU along the X and Y axis is high, the solution can  be potentially improved by selecting a different TU. TUs with a broader base can contain more boxes, potentially reducing the total number of TUs used. On the other side, TUs with a narrow base are used better and their cost is potentially lower.
\end{itemize}

After TUs destruction, algorithm 3DBP is repeatedly called with a different TU type taken from \TUS. Namely, \TUS\ is sorted by increasing volume $V_{tut}=X_{tut}\times Y_{tut}\times Z_{tut}$, then is circularly scanned with a pointer starting over at the first place as it passes the end of the list.
The pointer is initialized at the beginning of the list when Algorithm ILS-ATUCP is started and it determines the type of TU to use each time algorithm 3DBP is called (to build the first solution), during move $N3$ of LS1, and during LS2.
 
Local search LS2 stops when a call to 3DBP succeeds in finding an improving solution, in which case algorithm jumps back to LS1; in this case the pointer will determines the type of TU to pass to algorithm 3DBP during move $N3$, and the point to start from the next time LS2 will be executed. Otherwise, pointer completes a whole scan of \TUS\ and the algorithm definitely ends it run.

\subsection{Numerical example}\label{sec:num_example}

\begin{figure}[H]
	\caption{Boxes and EPs}
	\includegraphics[scale=0.60]{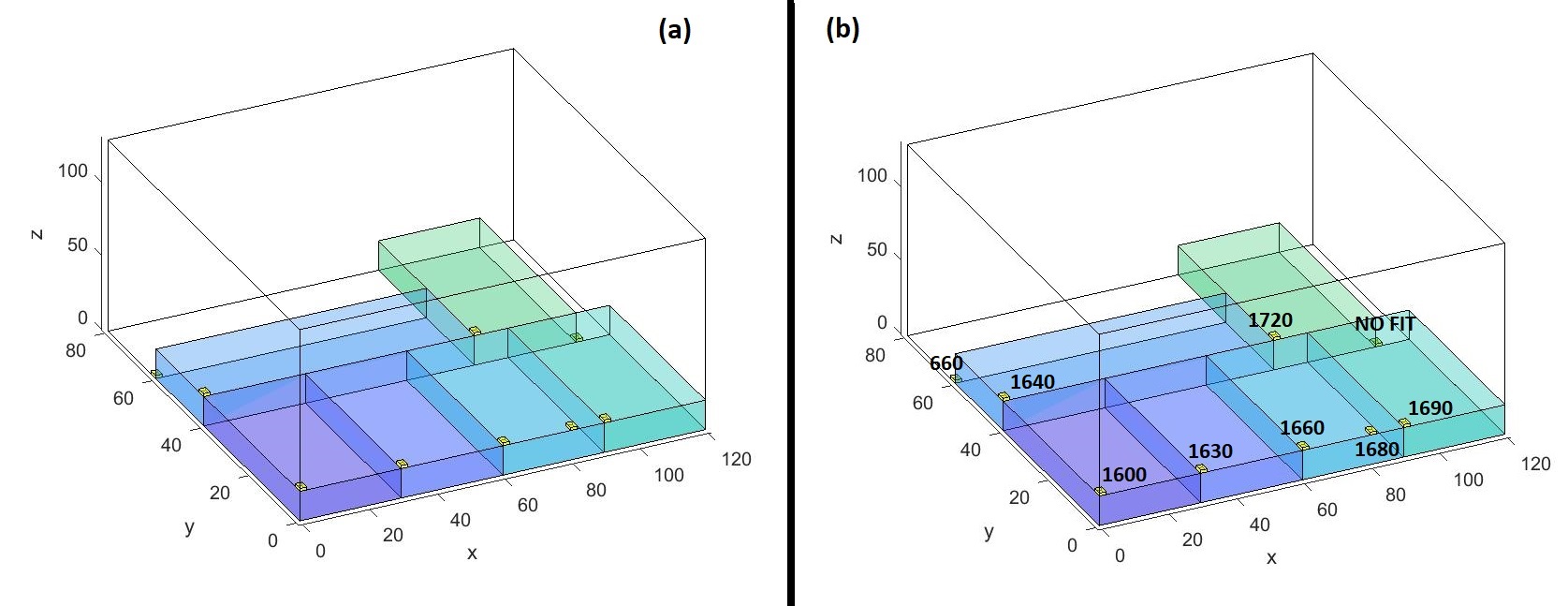}\label{fig:prescatola}
	\centering
\end{figure}

We present a  numerical example to better explain how the 3DBP algorithm works.
Figure \ref{fig:prescatola}(a) shows a partial solution produced by the 3DBP algorithm for a standard Europallet 120x80x130 where we highlight the currently available EPs. In the next iteration we need to insert a box of dimensions 30x40x20, weight is not considered in the example.

Table \ref{t:ep_esempio_misure} show the available EPs and the maximum dimensions allowed for a box to be placed at the EP.

 \begin{table}[H]
 	\centering
 	\caption{EPs coordinates and maximum measures} \label{t:ep_esempio_misure}
 	\resizebox*{0.30\textwidth}{!}{
 		\begin{tabular}[H]{l|c|c|c|c|c|c|}
  EP & $X_e$ & $Y_e$ & $Z_e$ & $MX_e$ & $MY_e$ & $MZ_e$  \\ \hline
1	&	30	&	0	&	20	&	90	&	80	&	110	\\ \hline
2	&	60	&	0	&	20	&	60	&	80	&	110	\\ \hline
3	&	0	&	60	&	0	&	80	&	20	&	130	\\ \hline
4	&	0	&	0	&	20	&	120	&	80	&	110	\\ \hline
5	&	90	&	0	&	20	&	30	&	80	&	110	\\ \hline
6	&	0	&	40	&	20	&	120	&	40	&	110	\\ \hline
7	&	80	&	0	&	20	&	40	&	80	&	110	\\ \hline
8	&	80	&	40	&	20	&	40	&	40	&	110	\\ \hline
9	&	110	&	40	&	0	&	10	&	40	&	130	\\

 		\end{tabular}
 	}
 \end{table}

It can be easily  checked that the box can fit in every EP with the only exception of EP 9, since the box has no dimension lower than 10, the maximum width allowed by the EP. For every EP where the box can fit, we calculate the cost of the EP using (\ref{eq:algo_obj}) for every available rotation and select the best rotation. We then select the EP associated with the minimum cost. Table \ref{t:ep_esempio_costi} shows the results of the calculation, represented visually by Figure \ref{fig:prescatola}(b). For simplicity, we only consider and report the cost of the best orientation for every EP.

 \begin{table}[H]
	\centering
	\caption{EPs coordinates and cost} \label{t:ep_esempio_costi}
	\resizebox*{0.30\textwidth}{!}{
		\begin{tabular}[H]{l|c|c|c|c}
EP & $X_e$ & $Y_e$ & $Z_e$ & Cost   \\ \hline
1	&	30	&	0	&	20	&	1630	\\ \hline
2	&	60	&	0	&	20	&	1660	\\ \hline
3	&	0	&	60	&	0	&	660	\\ \hline
4	&	0	&	0	&	20	&	1600	\\ \hline
5	&	90	&	0	&	20	&	1690	\\ \hline
6	&	0	&	40	&	20	&	1640	\\ \hline
7	&	80	&	0	&	20	&	1680	\\ \hline
8	&	80	&	40	&	20	&	1720	\\ \hline
9	&	110	&	40	&	0	&	NO FIT	\\
	
		\end{tabular}
	}
\end{table}

\begin{figure}[H]
	\caption{New solution once the box is inserted}
	\includegraphics[scale=0.35]{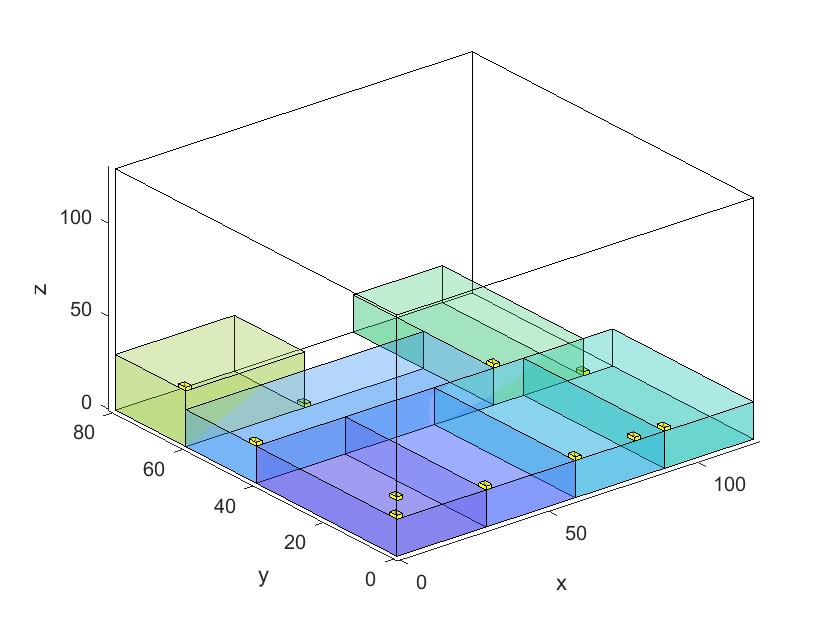}\label{fig:postscatola}
	\centering
\end{figure}

Our goal is to place the box in the lowest, western southern EP available. Even by visual inspection it is easy to verify that the best EP is (0,60,0) and the cost function reflects that. We note that the high difference between EP3 and all remaining EPs is due to the fact that EP3 is the only one with $z_e=0$. Figure \ref{fig:postscatola} shows the visual representation of the TU once we insert the box in the lowest-cost EP. In order to insert the box in EP (0,60,0) we must rotate it to have dimensions 40x20x30.

\section{Computational  results on 3D-BPP benchmark instances}\label{sec:results3DP}

In order to study the behaviour and the effectiveness of the ATUCP algorithm, the instances introduced by \citet{bischoff1995issues} for the 3D-BPP are solved. Seven subsets of instances do not consider weight of goods, while 5 subset of instances do. Each subset contain 100 instances and it is defined by the quantity of different types of boxes considered. For the sake of easing the comparison between the results, in this section the same output reported in \citet{bischoff1995issues} for the first 700 instances (the ones without considering weight) are compared with the results obtained by the 3DBP-Algorithm proposed in this paper. We now provide some details on the instances\\
There are seven subset of instances, each one containing 100 instances. In every instance boxes are regular hexahedrons. The seven subsets of instances are differentiated on the basis of the number of types of different boxes: 3, 5, 8, 10, 12, 15 and 20, respectively. For all instances a single type of bin is used with dimensions 587x233x220 (i.e. a 20 feet standard container). \\
The algorithm proposed in \citet{bischoff1995issues} is focused on maximizing the volume utilization of the container and ensuring the stability of boxes. It is also optimized to handle single bins and weakly heterogeneous boxes, while the algorithm presented in this work tries to find optimal solution using weakly heterogeneous bins and strongly heterogeneous boxes.\\
Table \ref{ATUCP_t:bench_results} show the results obtained for each subset of 100 instances solved in terms of: solution time (ST), average utilization rate of fullest container (A\%), highest utilization rate of the fullest container (Max\%), minimum utilization rate of the fullest container (Min\%).

\begin{table}[h!]
	\centering
	\caption{Computational results on 3D-BPP instances} \label{ATUCP_t:bench_results}
	\resizebox*{0.70\textwidth}{!}{
		\begin{tabular}[h!]{|l|c|c|c|c|}
			\hline
			Ist & ST(s) & A\% & Max\%  & Min\%  \\ \hline
			1&10.64&77.23\% (81.76\%)&90.01\% (94.36\%)&59.34\% (64.96\%)\\ \hline
			2&9.24&72.47\% (81.70\%)&83.12\% (93.76\%)&57.84\% (66.90\%)\\ \hline
			3&11.46&67.72\% (82.98\%)&79.99\% (92.76\%)0&56.96\% (66.91\%)\\ \hline
			4&10.44&69.98\% (82.60\%)&80.01\% (88.89\%)&59.32\% (66.46\%)\\ \hline
			5&13.12&69.01\% (82.76\%)&79.85\% (90.39\%)&61.01\% (70.38\%)\\ \hline
			6&11.75&67.85\% (81.50\%)&79.92\% (89.15\%) &57.44\% (64.86\%)\\ \hline
			7&12.89&67.77\% (80.51\%)&78.37\% (88.18\%)&62.12\% (70.50\%)\\ \hline
			
		\end{tabular}
	}
\end{table}

When compared to the results obtained in \citet{bischoff1995issues} reported in round parenthesis, the ATCUP algorithm performs reasonably, even if it is less competitive. The average utilization rate is lower and both maximum and minimum utilization rates are lower. The main reason for this performance issues are to be searched in the logic behind the 3D-BPP proposed and the problem tackled. In \citet{bischoff1995issues} the container loading problem is studied where the objective is to maximize the utilization rate of containers loaded, while in ATUCP the objective is different and considers multiple aspects of air transportation. Furthermore, the benchmark instances introduced by \citet{bischoff1995issues} are generated in order to study the maximum number of boxes that the algorithm can load in a single container. That means that the final solution  always contains at most two TUs, and therefore the only local search operators used in the ATUCP algorithm are the ones in $N1$ and $N2$. As we will see, the most effective local search operators are $N3$ and $LS2$. Thus, the improvements made by the use of local searches on the instances proposed in \citet{bischoff1995issues} are negligible. Overall, though, the ATUCP algorithm finds compact and good solutions. Utilization rate is satisfying and boxes are stable and few empty spaces are left between packed boxes.

\section{Experimental campaign}\label{sec:res}

No benchmark instances exist for the ATUCP. Therefore, we generated them according to the method explained in the following. We carried out two sets of experiments.

The first set aims at evaluating the solutions provided by  ILS-ATUCP  on instances for which a lower bound on the solution value can be determined. 
First we identify the lower bound solution for 100 instances (more details in the Appendix in Section \ref{sec:instances}). Then, starting from the lower bound solutions found, three partitioning schemes are used to obtain the list of boxes to be inserted in the TUs. The three algorithms partition the TUs of the lower bound solutions in order to generate the list of boxes boxes. The first scheme follows a layer based procedure to generate boxes. The second scheme generates boxes by filling up the volume of the TU. The third one generates a number of boxes of fixed dimensions that perfectly fit the TUs. Finally, the solution found by  ILS-ATUCP  is compared to the lower bound solution previously identified, in terms of number of TUs used and volume occupation.\\
The second set of tests aims at comparing the solution of ILS-ATUCP algorithm with the ones provided by the internally developed solver used by the Italian freight forwarding company which inspired this work. The comparison is made on two sets of instances. The first set consists of instances generated as described above. The second set consists of real instances provided by the company.

Before presenting the results, we now report some observations related to the parameters included in the objective function \eqref{for:OF}. According to the instances characteristics (described in the following) we set the value of $\Theta$ equal to 100: this provides solutions in which the value of the center of gravity term of the objective function is balanced with respect to the other two terms. The value of $\alpha$ has been set to 1 and we performed a sensitivity analysis at at verifying the robustness of solutions with respect to this parameters. As mentioned in Section \ref{sec:res_atucp_random}, this parameter has a negligible impact on the solution structure. Finally, the value of $\beta$ has been set to 100. Also in this case, when changing the value of $\beta$ the solutions structure did not change.

This section is organized as follows. Section \ref{sec:instances} describes how we generated ATUCP instances. Section \ref{sec:res_atucp_random} is devoted to the evaluation of ILS-ATUCP performance by comparison with the optimal solutions and to the sensitivity analysis on algorithm's parameters. Section \ref{sec:res_real} presents a comparison between algorithm's solutions and company's practice.

\subsection{Partitioning schemes to generate boxes}\label{s:part_alg}

Three partitioning schemes are proposed to partition the TUs identified by the optimal solution of the lower bound into boxes. 
They do not consider weight as the weight density is constant. Thus, they are just based on volumes.
\begin{itemize}
	\item The first scheme partitions the volume of each TU according to a layer-based procedure. The TU is first divided in layers based on the height, then each layer is divided first along one direction and then on the other. The sketch of the procedure is given in  Partitioning Scheme \ref{algo:part1}. $Z_{left},X_{left},Y_{left}$ are the residual height, width and length of the TU while $Z_{UB}, Y_{UB}, X_{UB}$ are the upper bound values of box height, length and width dimension respectively, $  Z_{LB}, Y_{LB}, X_{LB}$ are the lower bound values of box height, length and width dimension respectively. Upper and Lower bound values are inserted in order to limit the randomness of the generation procedure.
\begin{algorithm}[H]
	\SetAlgoLined
	Inizialization ($Z_{left} \gets Z_{TU}$) \\
	\While{$Z_{left} > 0$}{
		$X_{left} \gets X_{TU} $\\
		\eIf{$Z_{left} < Z_{LB}$}{
		$Z_b \gets Random (Z_{LB}; MIN(Z_{UB}; Z_{left}))$
		
	}{
$Z_b \gets Z_{left}$
}

	\While{$X_{left} > 0$}{
		$ Y_{left} \gets Y_{TU}$\\
	\eIf{$X_{left} < X_{LB}$}{
		$X_b \gets Random (X_{LB}; MIN(X_{UB}, X_{left})$
		
	}{
		$X_b \gets Z_{left}$
	}

	\While{$Y_{left} > 0$}{
	\eIf{$Y_{left} < Y_{LB}$}{
		$Y_b \gets Random (Y_{LB}; MIN(Y_{UB}, Y_{left}))$\\
		Generate box with dimensions ($X_b, Y_b, Z_b$)\\
		
	}{
		$Y_b \gets Y_{left}$ \\
		Generate box with dimensions ($X_b, Y_b, Z_b$)\\
	}
	update $Y_{left}$\\	
}
update $X_{left}$\\	
}
		
	update $Z_{left}$\\		
	
	}
 \caption{Partitioning Scheme 1: Layer-based procedure}\label{algo:part1}
\end{algorithm}
Figure \ref{fig:layer} shows a graphical illustration of the partitioning scheme 1, where layer 1 is the layer identified along the Z axis, layer 2 the layer  along the X axis and Layer 3 is the last one,  along the Y axis. 
 \begin{figure}[H]
	\caption{Partitioning Scheme 1: Graphical representation}
	\includegraphics[scale=0.35]{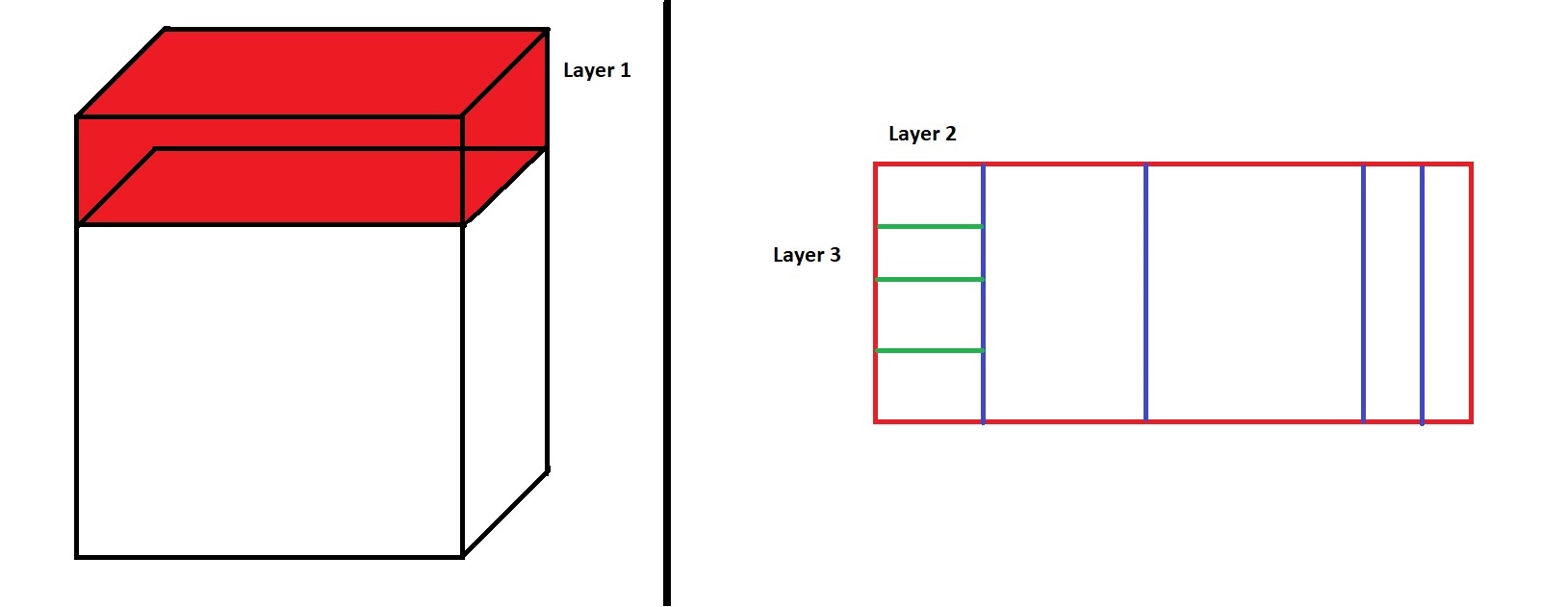}\label{fig:layer}
	\centering
\end{figure}

	\item The second scheme starts by generating a box with random dimensions in the southern-eastern-down corner of the TU, the same location where EP (0,0,0) is positioned when the TU is empty. Once the box is generated, one of the faces directed towards the inside of the TU is selected at random (direction 1, 2 or 3 in Figure \ref{fig:scava}). Then, boxes are generated  with the dimension of the face selected fixed, and random dimension along the direction faced by the selected face of the box. The partitioning scheme partitions the "tunnel" generated by projecting the selected face in the direction of the projection. For instance, direction 2 projects the face in the YZ plane along the X axis, so it partitions the projected tunnel along the X axis. The boxes generate will have identical $Y_b$ and $Z_b$ dimensions and different $X_b$ dimension. Once one direction is fully partitioned, new EPs are generated and the procedure continues until the volume of the TU is fully partitioned. %Once a box is generated, the weight is assigned based on the volume and total weight of the instance.
	The procedure is described in Partitioning Scheme \ref{algo:part2} where $f$ is the selected face.

	\begin{algorithm}[H]
		\SetAlgoLined
		Inizialization: Add EP (0,0,0) to EP-List \\
		\While{$|EP-List| >0$}{
select one EP $e$ at random from EP-List;\\
Generate one box $b$ with random dimensions ($X_b, Y_b, Z_b$);\\
$f \gets Random(1,3)$;\\
\If{f = 1}{
	\While{$Z_{left} > 0 $}{
		$Z_u \gets Random (Z_{LB}, Z_{left})$;\\
Generate one box $u$ with dimensions ($X_b, Y_b, Z_u  $)\\
Update $Z_{left}$
}
Remove $e$ from EP-List, add EPs in position 1 and 2 (if available) to EP-list according to Figure \ref{fig:scava}\\
}

\If{f = 2}{
		\While{$X_{left} > 0 $}{
		$X_u \gets Random (X_{LB}, X_{left})$;\\
		Generate one box $u$ with dimensions ($X_u, Y_b, Z_b  $)\\
		Update $X_{left}$}
	Remove $e$ from EP-List, add EPs in position 1 and 3 (if available) to EP-list according to Figure \ref{fig:scava}\\
}
\If{f = 3}{
		\While{$Y_{left} > 0 $}{
		$Y_u \gets Random (Y_{LB}, Y_{left})$;\\
		Generate one box $u$ with dimensions ($X_b, Y_u, Z_b  $)\\
		Update $Y_{left}$}
	Remove EP $e$ from EP-List, add EPs in position 2 and 3 (if available) to EP-list according to Figure \ref{fig:scava}\\
}
		}
		\caption{Partitioning Scheme 2}\label{algo:part2}
	\end{algorithm}

Figure \ref{fig:scava} shows a visual representation of the idea behind the second partitioning scheme. Once a box is created, one direction is chosen and the "tunnel" created by projecting the face directed towards the chosen direction is partitioned into boxes until the edge of the TU is reached. From there, new EPs are added and used as starting location for the generation of new boxes and the procedure iterates until the volume of the TU is fully partitioned. %The weight is assigned to the boxes in order to maintain the original weight density.

 \begin{figure}[H]
	\caption{Partitioning scheme 2: Graphical representation}
	\includegraphics[scale=0.35]{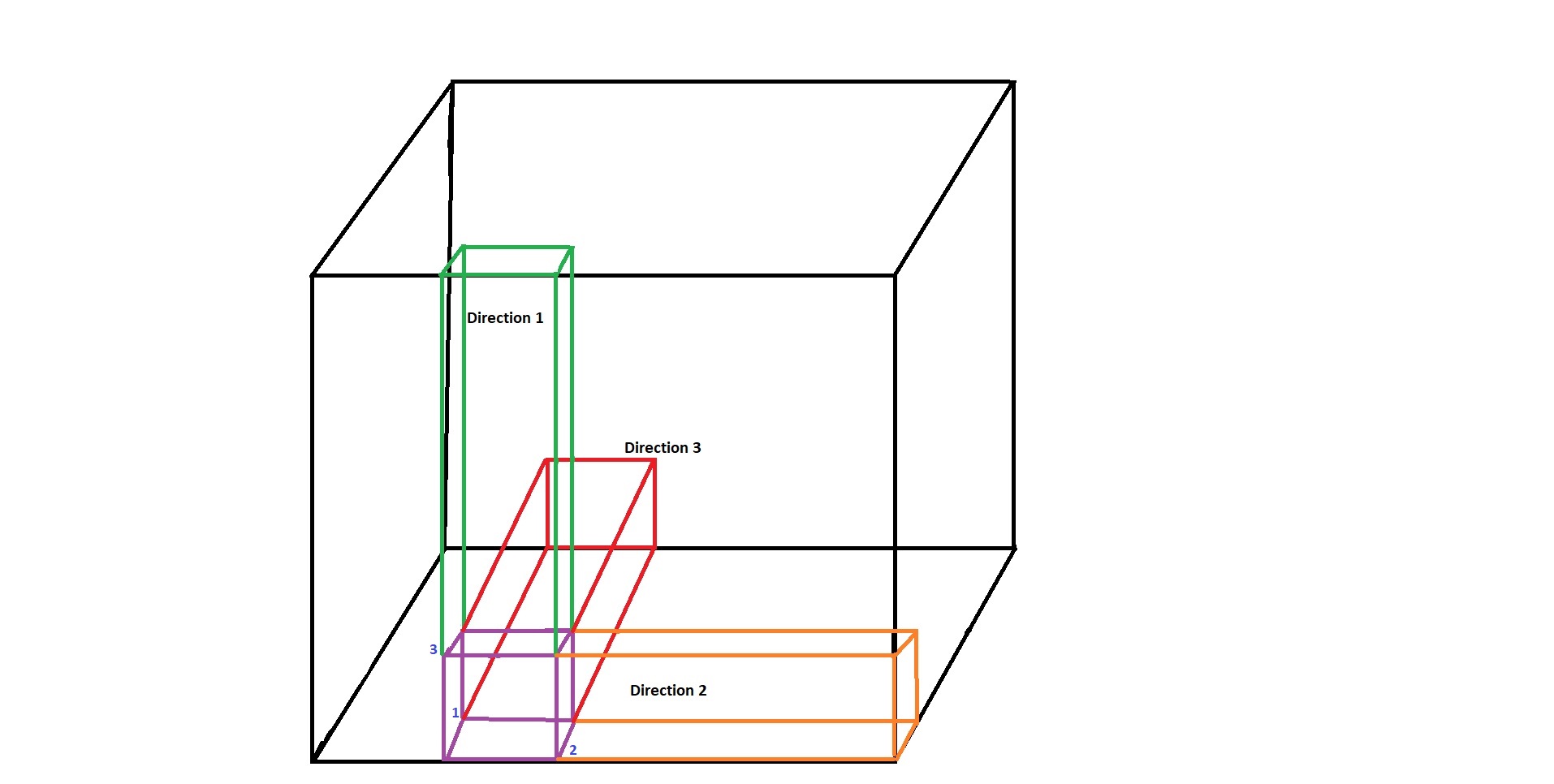}\label{fig:scava}
	\centering
\end{figure}

	\item The third partitioning scheme identifies a perfect partition of each type of TU using boxes of standard dimensions and simply generates boxes based on the lower bound solution. For instance, a pallet 120x80x160 can be partitioned in 32 boxes 40x30x40, so for each TU with dimension 120x80x160 in the lower bound solution, 32 boxes 40x30x40 are generated. The dimensions of the boxes are chosen arbitrarily to perfectly partition the volume of the TU. The number and dimension of the boxes generated is the same for each kind of TU used in the solution, i.e. a TU 120x80x130 will generate a different number of boxes with different dimensions than those generated by a TU 120x100x160, but every TU 120x80x130 will generate the same boxes.  %The weight is assigned to each box respecting the weight density of the instance.
	
\end{itemize}

Table \ref{t:part_data} in the Appendix reports some statistics about the instances generated by the three partitioning schemes.
The first scheme generates the lowest number of boxes, but the variety of boxes dimensions is the largest, both in terms of number of unique dimensions and medium square deviation from the average for each dimension (the average square deviations for the three dimensions are (30.21, 27.65, 32.82) for partitioning scheme 1, (11.42, 9.54, 7.10)  for partitioning scheme two and (0, 3.15, 1.07) for partitioning scheme three). The second scheme generates a high number of small boxes with a great variety of dimensions. The third scheme  generates an high number of boxes, often higher than the second partitioning scheme, but the number of unique dimensions is limited and the dimensions are more suitable for simpler and clearer allocations of the boxes in the TUs.

\subsection{Computational results on randomly generated instances}\label{sec:res_atucp_random}

The 300 instances generated (100 instances for each partitioning scheme) are solved with ILS-ATUCP  and the solutions are  compared with the lower bounds.
Moreover, different values of $\Omega$ in LS2 are used to solve the 300 instances and the solutions are studied in order to carry out a sensitivity analysis. We also performed experiments by varying the value of $\gamma$ in LS2 but noticed that this parameter has no impact on the solution produced by LS2. Thus, we fixed it to 100 in the following experiments.
In Table \ref{t:risultati_ATUPC} some aggregated statistics statistics for the 100 instances generated using each partition algorithm are reported.

\begin{table}[H]
	\centering
	\caption{Exact methodology results} \label{t:risultati_ATUPC}
	\resizebox*{0.35\textwidth}{!}{
		\begin{tabular}[h!]{|l|l|c|c|c|}
			\hline
		$\Omega$&Stat&A1&A2&A3\\ \hline
\multirow{9}{*}{75}&
N° TU	&	12.06	&	14.09	&	15 \\
&$\Delta$ Volume 	&	11.06 \%	&	31.08\%	&	15.52\% \\
&$CI^{xy}$	&	0.114	&	0.119	&	0.002 \\
&$CI^z$ & 0.485 & 0.467 & 0.499 \\
&Sol. Time (s)	&	1.64	&	157.69	&	2.04 \\
&N° Opt	&	28	&	6	&	4 \\
&MAX \%	&	87.29 	&	83.88	&	92.92 \\
&MIN \%	&	70.53	&	71.32	&	88.12 \\
%&MAX \%	&	100	&	100	&	100 \\
 \hline

 \multirow{9}{*}{80}&
N° TU	&	10.63	&	11.39	&	15 \\
&$\Delta$ Volume 	&	14.57 \%	&	35.64\%	&	15.52\% \\
&$CI^{xy}$	&	0.113	&	0.117	&	0.002 \\
&$CI^z$ & 0.499 & 0.575 & 0.501  \\
&Sol. Time	&	1.98	&	188.57	&	2.04	\\
&N° Opt	&	24	&	5	&	5 \\
&Avg. MAX \%	&	88.17 	&	84.12	&	92.92 \\
&Avg. MIN \%	&	70.55	&	71.48	&	88.12 \\
%&MAX \%	&	100	&	100	&	100 \\
 \hline

 \multirow{9}{*}{85}&
N° TU	&	9.69	&	11.13	&	15 \\
&$\Delta$ Volume 	&	12.48 \%	&	34.42\%	&	15.52\% \\
&$CI^{xy}$	&	0.114	&	0.109	&	0.002 \\
&$CI^z$ & 0.481  & 0.466 & 0.501 \\
&Sol. Time	&	2.47	&	303.13	&	2.04	\\
&N° Opt	&	25	&	7	&	5 \\
&MAX \%	&	88.06 	&	84.27	&	92.92 \\
&MIN \%	&	70.22	&	72.12	&	88.12 \\
%&MAX \%	&	100	&	100	&	100 \\
 \hline

 \multirow{9}{*}{90}&
N° TU	&	9.08	&	11.48	&	12.98 \\
&$\Delta$ Volume 	&	12.10 \%	&	34.02\%	&	10.01\% \\
&$CI^{xy}$	&	0.115	&	0.111	&	0.002 \\
&$CI^z$ & 0.488 & 0.475 & 0.501 \\
&Sol. Time	&	3.03	&	303.38	&	2.04	\\
&N° Opt	&	23	&	7	&	6 \\
&MAX \%	&	88.09	&	84.43	&	93.87 \\
&MIN \%	&	71.11	&	72.25	&	88.84 \\
%&MAX \%	&	100	&	100	&	100 \\
 \hline

 \multirow{9}{*}{95}&
N° TU	&	9.07	&	11.35	&	10.06 \\
&$\Delta$ Volume 	&	12.09 \%	&	33.94\%	&	4.22\% \\
&$CI^{xy}$	&	0.116	&	0.109	&	0.001 \\
&$CI^z$ & 0.495 & 0.498 & 0.501 \\
&Sol. Time	&	2.76	&	302.54	&	11.53	\\
&N° Opt	&	24	&	7	&	21 \\
&MAX \%	&	88.35 	&	84.55	&	99.77 \\
&MIN \%	&	71.34	&	72.31	&	91.78 \\
%&MAX \%	&	100	&	100	&	100 \\
\hline

		\end{tabular}
	}
\end{table}

Columns A1, A2 and  A3 refers to the set of instances generated by the first, second and third partitioning scheme, respectively.
The first row reports the average number of TUs used in the solutions found, the second row is the average difference of volume generated with respect to the lower bound solutions. For difficult instances the ILS-ATUCP algorithm is incapable of perfectly recreating the original TUs layout, the solution found generates a higher volume.
The third and fourth row reports the average value of the center of gravity coordinates $CI^{xy}$ and $CI^z$ as determined in \eqref{for:OF}, respectively.
The fifth row reports the average solution time in seconds. The sixth row reports the number of instances where the solution found by the ILS-ATUPC algorithm is equal to the lower bound solutions. The seventh and eighth row report the average maximum and minimum filling percentage of the solutions found, respectively.

It is interesting to note the evolution of the solution when varying the value of $\Omega$.
For the first and second partitioning schemes, columns A1 and A2, the pattern is fairly clear. By allowing a wider search in LS2, the ILS-ATUCP algorithm is capable of reducing the number of TUs used. On the other side, the volume developed is slightly higher on average. This is because while noticeable improvements are made by the use of less, bigger TUs, in some instances the search deviates highly and find worse solutions, since some TUs with high fill rate (85+\%) are destroyed and rebuilt with a different kind of TU which is not always optimal. 
Moreover, this result is influenced by the value of the parameter $\beta$. The higher is this value, the greatest are the efforts  in reducing the number of TUs used, even if this leads to an increase in the volume used. The solutions of the instances generated using the third partitioning scheme present a particular behaviour. For every value of $\Omega$ with the exception of $\Omega =90$ and  $\Omega =95$, the solution found is always the same.
This is due to the fact that the solution found before LS2 is good enough to never trigger LS2, and so the solution space is not thoroughly explored. With $\Omega = 90$ and $\Omega = 95$, or by forcing the search through LS2, different dimensions of TUs are used and the solutions found show great improvements. The reduction in the number of TUs used and the difference in the gap with the solutions found by the exact method are measures of the improvements found through the use of LS2. Moreover, the center of gravity position found by the ILS-ATUPC algorithm for the instances generated by the third partitioning scheme is optimal in almost every TU built.

We now move to the analysis of the ILS-ATUCP algorithm performance on the basis of the partitioning scheme used to generate the instances.
The solutions found in the instances generated by the first partitioning scheme are such that the TUs used are compact and the volume is well covered. This is due to the layer structure of the instance. 
While it is nearly impossible for the algorithm to perfectly recreate the original layers, it proves its effectiveness in creating TUs that are compact, with boxes packed together in order to create balanced, leveled layers with few empty spaces between boxes and towards the edge of the TU. The center of gravity is in satisfying position in almost every solution, both in terms of base and height position.

\begin{figure}[H]
	\caption{Example of solutions found for an instance generated using Partitioning Scheme 3}
	\includegraphics[scale=0.45]{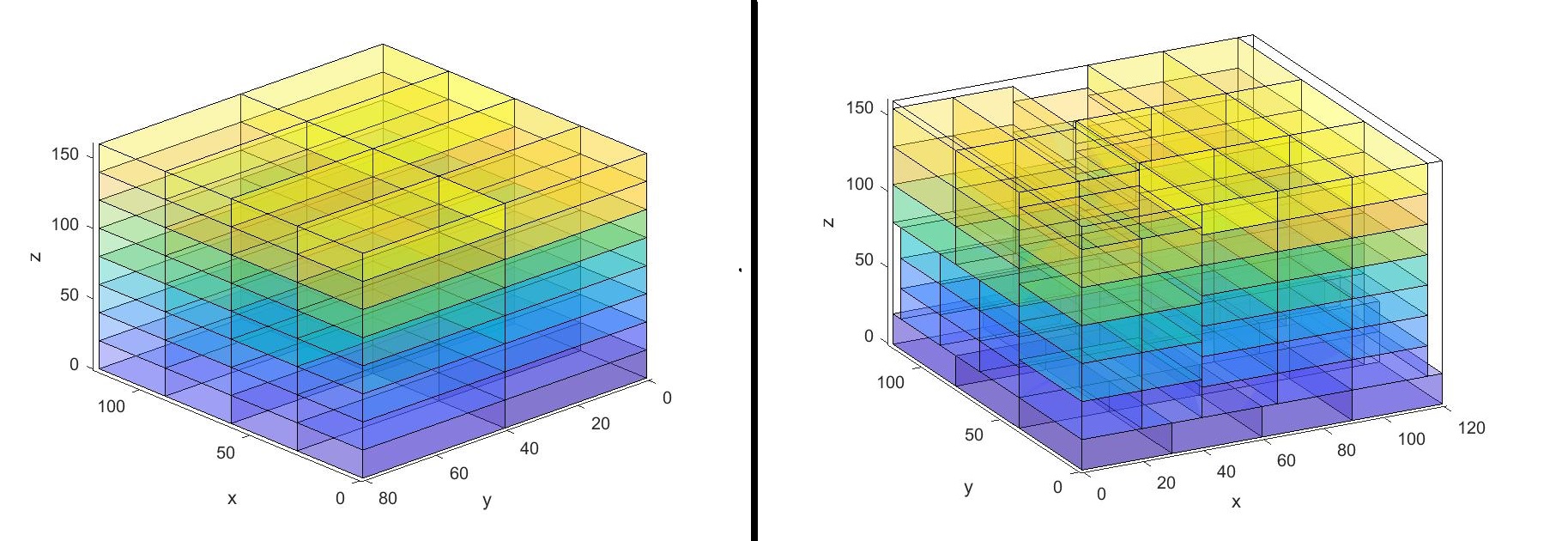}\label{fig:Bancalozzi}
	\centering
\end{figure}

The second partitioning scheme, on the other side, produces a high number of boxes of almost unique dimensions each, and building good TUs proves to be a challenging task, especially on instances with high weight and volume. 
Considering the difficulty of the task, the solutions generated by the ILS-ATUPC algorithm are compact and with few empty spaces. Moreover, the position of the center of gravity is more centered and at a lower height if compared with the instances generated with the first partitioning scheme. The time required to find a solution, on average, is considerably higher and the impact of widening the research by increasing the value of $\Omega$ is evident. On the other side, the improvements in the solutions found, made possible by a further use of LS2, shows the effectiveness of the local search.\\
Solutions of the instances generated by the third partitioning scheme are a good example of the potential effectiveness in real life scenarios. Solutions are found in fast times and the TUs built are compact, balanced and with a high utilization rate of the TU volume. Figure \ref{fig:Bancalozzi} shows and example of two TUs built on an instance generated using the third partitioning scheme. One TU is perfectly reproduced by using boxes of equal dimensions that fits the TU shape, while the other contain boxes of various dimensions and are assembled together in a compact, stable and centered layout.

\begin{figure}[H]
	\caption{Example a solution found for an instance generated by Partitioning Scheme 1}
	\includegraphics[scale=0.30]{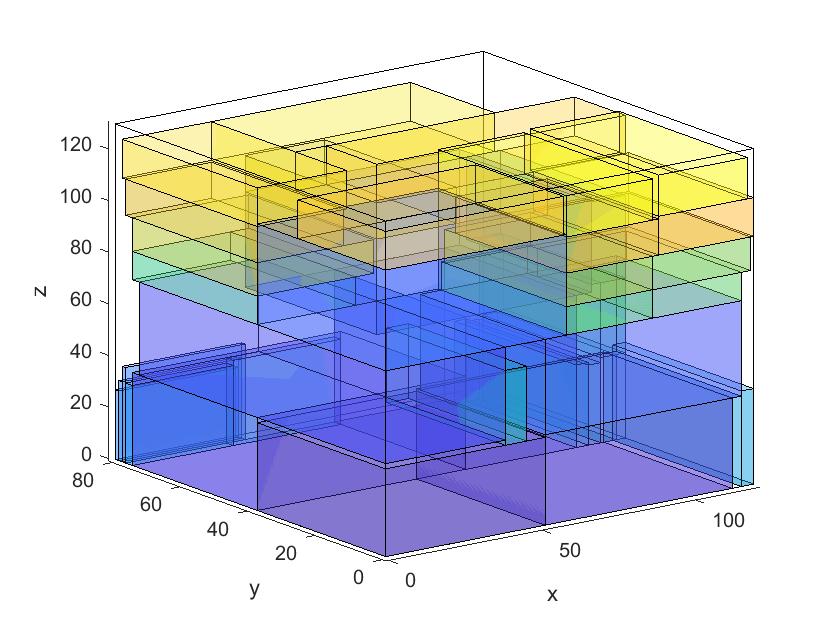}\label{fig:banc_diff}
	\centering
\end{figure}

\begin{figure}[H]
	\caption{Example of solutions found for an instance generated by Partitioning Scheme 2}
	\includegraphics[scale=0.45]{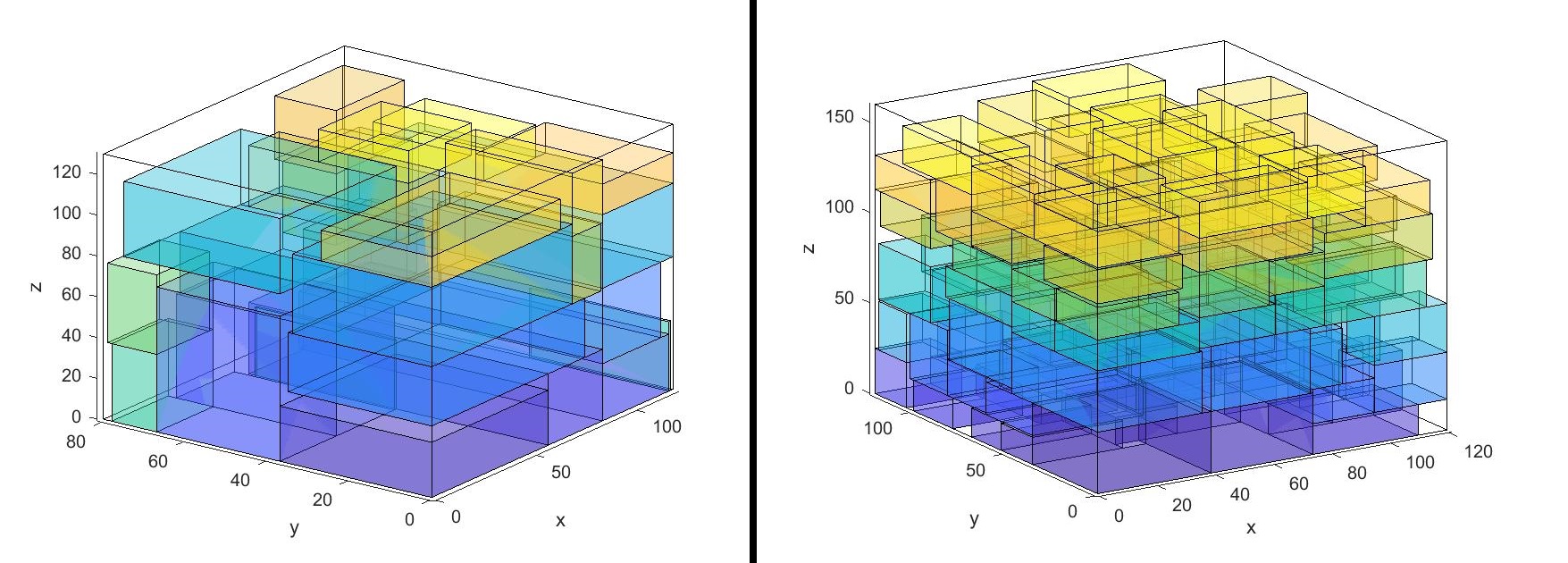}\label{fig:scavarisultati}
	\centering
\end{figure}

Figure \ref{fig:banc_diff} and \ref{fig:scavarisultati} show a TU built from an instance generated with the first and second partitioning scheme, respectively. Even though the number of boxes is high and each of them has unique dimensions, the ILS-ATUPC algorithm finds a good solution that nicely occupy the volume of the TU.

Overall, the TUs built are compact, stable, found in relatively fast times and have a satisfying utilization rate.

The effectiveness of the local search procedures is reported in Table \ref{t:local_results}. The first row reports the number of instances where LS1 finds improvements on the original solution. The second and third rows report the average and max improvements found by LS1. The fourth row report the number of instances where LS2 finds improvements of the incumbent solution. The fifth and sixth rows report the average and maximum improvement found by LS2. The last two rows are the average and max reduction in the number of TUs obtained through both local search procedures. As expected, LS1 never reduces the number of TUs, while LS2 is fairly effective in that. LS1 finds little upgrades, mostly through N1 and N3, while N2 is difficult to apply due to the fact that is difficult to find boxes that can be swapped without breaking the overlapping constraint. LS2 on the other side respects the expectations and is fairly efficient in improving the solution.

\begin{table}[H]
	\centering
	\caption{Local search procedures analysis} \label{t:local_results}
	\resizebox*{0.35\textwidth}{!}{
		\begin{tabular}[h!]{|l|l|c|c|c|}
			\hline
			$\Omega$&Stat&A1&A2&A3\\ \hline
			\multirow{8}{*}{75}&
			N° Improvements LS1	&	57	&	55	&	0\\
			& Avg. Improvement LS1 \%	&	1.23 	&	1.45	&	0 \\
			& MAX improvement LS1 \%	&	1.78	&	2.01	&	0\\
			&N° Improvements LS2	&	94	&	98	&	0 \\
			&Avg. Improvement LS2 \% &	9.67 	&	9.97 	&	0 \\
			& MAX improvement LS2 \%	&	41.22 	&	33.44	&	0 \\
			&Avg TU's Reduction \%	&	17.23 	&	16.73	&	0 \\
			&MAX TU's reduction \%	&	50	&	50	&	0\\
			
			\hline
			
			\multirow{8}{*}{80}&
			N° Improvements LS1	&	57	&	55	&	0\\
& Avg. Improvement LS1 \%	&	1.23 	&	1.45	&	0 \\
& MAX improvement LS1 \%	&	1.78	&	2.01	&	0\\
			&N° Improvements LS2	&	94	&	98	&	0 \\
			&Avg. Improvement LS2 \% &	9.72 	&	9.87 	&	0 \\
			& MAX improvement LS2 \%	&	41.22 	&	33.44	&	0 \\
			&Avg TU's Reduction \%	&	17.23 	&	16.73	&	0 \\
			&MAX TU's reduction \%	&	50	&	50	&	0\\

\hline
			\multirow{8}{*}{85}&
			N° Improvements LS1	&	57	&	55	&	0\\
& Avg. Improvement LS1 \%	&	1.23 	&	1.45	&	0 \\
& MAX improvement LS1 \%	&	1.78	&	2.01	&	0\\
			&N° Improvements LS2	&	96	&	99	&	0 \\
			&Avg. Improvement LS2 \% &	9.73	& 9.98	&	0 \\
			& MAX improvement LS2 \%	&	44.54	&	35.15	&	0 \\
			&Avg TU's Reduction \%	&	17.74 	&	16.98	&	0 \\
			&MAX TU's reduction \%	&	50	&	50	&	0\\

\hline
			
			\multirow{8}{*}{90}&
			N° Improvements LS1	&	57	&	55	&	0\\
& Avg. Improvement LS1 \%	&	1.23 	&	1.45	&	0 \\
& MAX improvement LS1 \%	&	1.78	&	2.01	&	0\\
			&N° Improvements LS2	&	97	&	99	&	35 \\
			&Avg. Improvement LS2 \% &	9.76	& 9.78	&	11.21 \\
			& MAX improvement LS2 \%	&	42.32	&	33.79	&	15.63 \\
			&Avg TU's Reduction \%	&	17.32 	&	16.81	&	9.12\% \\
			&MAX TU's reduction \%	&	50	&	50	&	50\\

\hline
			
			\multirow{8}{*}{95}&
			N° Improvements LS1	&	57	&	55	&	0\\
& Avg. Improvement LS1 \%	&	1.23 	&	1.45	&	0 \\
& MAX improvement LS1 \%	&	1.78	&	2.01	&	0\\
			&N° Improvements LS2	&	97	&	99	&	98 \\
			&Avg. Improvement LS2 \% &	9.79	& 9.98	&	18.02 \\
			& MAX improvement LS2 \%	&	41.54	&	33.15	&	35.53 \\
			&Avg TU's Reduction \%	&	17.34 	&	16.88	&	33.92 \\
			&MAX TU's reduction \%	&	50	&	50	&	50\\

\hline

		\end{tabular}
	}
\end{table}

Regarding the position of the center of gravity, the solutions found are satisfying. The average center of gravity coordinates are always close to the center of the TU for each axis, for each type of TU used. Indeed, the maximum deviation of the center of gravity coordinates with respect to the center of the TU is 2\% on the X-axis and 3.6\% on the Y-axis. For the Z axis the center of gravity depends on the TU occupation, but still the maximum deviation with respect to the half-height of the TU is 10.7\%. We performed a sensitivity analysis to evaluate whether the value of $\alpha$ has an impact on the structure of the solution, in terms of number and type of TUs used and in terms of positioning of boxes in the TUs, and thus center of gravity. In particular, we tested values of $\alpha =50,100,1000$. We noticed that solutions did not change, with respect to the ones found with $\alpha=1$, in terms of number and types of TUs used. Deviations of center of gravity positions were also negligible, lower than 0.1\% on all coordinates and all types of TUs, for all values of $\alpha$. Thus, we concluded that the ILS-ATUCP algorithm always provides stable solutions in terms of boxes disposition.

Table \ref{t:CG_pos} reports the average, minimum and maximum position for each coordinates of the center of gravity, throughout all the 300 instances solved, for different values of $\alpha$, i.e, the weight of the center of gravity components in function \eqref{for:OF}. The $CI^{xy}$ and $CI^z$ indexes present satisfying results, with values close to those of the theoretical optimum, which corresponds to the mid $X$, $Y$ and $Z$ values of the TU reported in the first column of the table.

\begin{table}[H]
	\centering
	\caption{Average center of gravity coordinates} \label{t:CG_pos}
	\resizebox*{0.65\textwidth}{!}{
		
		\begin{tabular}[h!]{|l|l|ccc|ccc|ccc|}
			\hline
			& & \multicolumn{3}{c}{X coordinates} &\multicolumn{3}{c}{Y coordinates} &\multicolumn{3}{c}{Z coordinates} \\ \hline
		$\alpha$ value	& TU & Avg. & Min & Max & Avg. & Min & Max & Avg. & Min & Max \\ \hline
			\multirow{6}{*}{$\alpha=1$}&
			120x80x130& 58.83 & 55.02 & 64.75 & 38.56 & 37.42 & 43.92& 60.44 & 48.12 & 75.25 \\
		&	120x80x160 & 58.77 & 55.01 & 64.73 & 38.57 & 37.55 & 43.89 & 76.12 & 70.98 & 98.44 \\
		&	120x100x130 & 58.82 &55.09 &65.01 &48.64 &44.18 & 55.32 & 59.21 & 48.06 & 73.43 \\
		&	120x100x160 & 58.75 &54.98 &64.73 &48.58 & 45.76 & 54.99 & 73.32 & 68.12 & 89.68\\
		&	120x120x130 &58.83 &55.03 &65.08 &58.67 & 55.33 &65.12 &58.02 & 47.29 &71.73\\
		&	120x120x160 & 58.85 &54.99 &64.12 &58.61 &56.01 & 66.23 &71.98 & 65.97 & 86.84\\ \hline
			
			\multirow{6}{*}{$\alpha=50$}&
            120x80x130& 58.83 & 55.02& 64.75 & 38.56 & 37.42 & 43.92& 60.44 & 48.12 & 75.25 \\
        &	120x80x160 & 58.77 & 55.01  & 64.73 & 38.57 & 37.55 & 43.89 & 76.12 & 70.98 & 98.44 \\
        &	120x100x130 & 58.82 &55.09 &65.01 &48.64 &44.18 & 55.32 & 59.21 & 48.06 & 73.43 \\
        &	120x100x160 & 58.75 &54.98 &64.73 &48.58 & 45.76 & 54.99 & 73.32 & 68.12 & 89.68\\
        &	120x120x130 &58.83 &55.03 &65.08 &58.67 & 55.33 &65.12 &58.02 & 47.29 &71.73\\
        &	120x120x160 & 58.85 &54.99 &64.12 &58.61 &56.01 & 66.23 &71.98 & 65.97 & 86.84\\ \hline

			\multirow{6}{*}{$\alpha=100$}&
            120x80x130& 58.57 & 55.22 & 64.55 & 38.76 & 37.55 & 43.72& 60.14 & 48.37 & 75.03\\
        &	120x80x160 & 58.68 & 55.19 & 64.59 & 38.74 & 37.63 & 43.71 & 75.97 & 71.16 & 98.23 \\
        &	120x100x130 & 58.77 &55.31 &64.97 &48.85 &44.36 & 55.22 & 58.94 & 48.32& 73.22 \\
        &	120x100x160 & 58.62 &55.09&64.54 &48.77 & 45.96 & 54.75 & 73.09 & 68.34 & 89.51\\
        &	120x120x130 &58.75 &55.17 &64.97 &58.91 & 55.55 &64.98 &57.83 & 47.55 &71.50\\
        &	120x120x160 & 58.75 &55.12 &64.00 &58.89 &56.22 & 66.05 &71.77 & 66.19 & 86.61\\ \hline

			\multirow{6}{*}{$\alpha=1000$}&
            120x80x130& 58.41 & 55.44 & 64.35 & 38.91 & 37.42 & 43.54& 59.99 & 48.61 & 74.93 \\
        &	120x80x160 & 58.37 & 55.47 & 64.37 & 38.92 & 37.55 & 43.49 & 75.85 & 71.32 & 98.01 \\
        &	120x100x130 & 58.46 &55.62 &64.84 &48.98 &44.18 & 55.02 & 58.79& 48.63 & 73.01 \\
        &	120x100x160 & 58.55 &55.31 &64.34 &48.97 & 45.76 & 54.61 & 72.98 & 68.68 & 89.27\\
        &	120x120x130 &58.39 &55.51 &64.81 &59.09 & 55.33 &64.86 &57.69 & 47.75 &71.23\\
        &	120x120x160 & 58.43 &55.34 &63.89 &59.04 &56.01 & 65.97 &71.65 & 66.56 & 86.29\\ \hline

		\end{tabular}
	}
\end{table}

The most noticeable aspect is that throughout all instances solved, the center of gravity is fairly centered along the X and Y axis. The height of the CG is on average slightly below half the maximum height allowed, which is fairly satisfying. One interesting fact is that by enlarging the base of the TU, the height of the center of gravity decreases. This means that by allowing a wider space on the base, the algorithm is capable of finding improved layouts with centers of gravity placed at a smaller height.\\
We carried out a sensitivity analysis by varying the value of parameter $\alpha$ in equation \eqref{for:OF}. The aim is to check whether, and how, the solution may vary by changing the weight given to the position of the center of gravity in the computation of the objective function. The values of $\alpha$ are set 1, 50, 100 and 1000. The results show the variation of the value of $\alpha$ brings minimal changes to the solutions found, meaning that the ATUCP algorithm always provides stable solutions in terms of boxes disposition.\\
This is mainly due to the fact that the search of the optimal placement of boxes in the TU, ensuring stable, compact layouts with centered centers of gravity, is given by the cost function \eqref{eq:algo_obj} of the 3DBP-Algorithm. By changing the parameters in cost function \eqref{eq:algo_obj}, different objectives are explored (for instance by changing the value of parameter $M$, the TUs generated have a center of gravity in lower position, but less centered on the base), but the solutions found are maintained substantially unvaried with the variation of parameter $\alpha$ of equation \eqref{for:OF}.

\subsection{Comparison with company's practice}\label{sec:res_real}

The second set of tests is about comparing the solution found by ILS-ATUCP  with the solution found by operators of an Italian freight forwarding company. Both real world and randomly generated instances were used for comparison.

Randomly generated instances are a selection of 10 instances generated with the partitioning scheme 1, 10 instances generated with the partitioning scheme 2, and five instances  generated with the partitioning scheme 3. The last ones have a low variety of types of boxes and dimensions are fairly standard and fit the TUs. They are used as a benchmark to verify whether the solutions provided by ILS-ATUCP  fits with the ones produced by company's operators. 

Results are shown in Table \ref{t:finte}. The first two columns report the number of TUs used by the company's operators and by ILS-ATUCP algorithm, respectively, while the last column reports the difference in the total volume generated by the ILS-ATUCP solution with respect to the solution of the company's operators. The last column reports the CPU time of ILS-ATUCP.

In the simpler instances (rows 1-5 in Table \ref{t:finte}), i.e., those generated with partitioning scheme 3, the performance of the ILS-ATUCP algorithm and operators are similar, thus proving that ILS-ATUCP is capable of producing solutions fitting with practical rules in a short computing time. 

Results for the 10 instances generated by the first partitioning scheme are reported in rows 6-15 of Table \ref{t:finte} and show that ILS-ATUCP largely outperforms company's practice, both in terms of number of TUs used and in terms of volume generated. Moreover, the algorithm is capable of finding the solution in very limited time.
Results for the instances generated with the second partitioning scheme are reported in rows 16-25. The algorithm again finds much better solutions in every instance. 

\begin{table}[H]
	\centering
	\caption{Comparison on randomly generated instances} \label{t:finte}
	\resizebox*{0.35\textwidth}{!}{
		\begin{tabular}[h!]{|l|c|c|c|c|c|}
			\hline
			I& N boxes &  TU Op & TU ILS-ATUCP & Vol $\Delta$ (\%) & Time ILS-ATUCP (sec.) \\ \hline
1	& 64 &	4	&	4	&	0 &	2\\
2	& 48 &	3	&	3	&	0 &	3\\
3	& 32 &	2	&	2	&	0 &	2\\
4	& 128 &	6	&	6	&	0 &	1\\
5	& 144 &	7	&	7	&	0 &	2\\
Avg. & 83.20 & 4.40& 4.40 & 0 & 2\\ \hline
%6	&	8	&	8	&	-3.12	\\
%7	&	7	&	7	&	-2.23	\\
%8	&	5	&	4	&	-6.43	\\
%9	&	4	&	4	&	1.22	\\
%10	&	5	&	6	&	3.12	\\
%11	&	10	&	9	&	-7.32	\\
%12	&	9	&	8	&	-9.83	\\
%13	&	6	&	7	&	1.32	\\
%14	&	3	&	4	&	1.12	\\
%15	&	20	&	18	&	-7.12	\\ \hline
6	& 233 &	10	&	8	& 	-9.12 & 2	\\
7	& 135 &	9	&	7	& 	-7.32 &	3\\
8	& 99 &	5	&	4	& 	-12.27 & 5	\\
9	& 76 &	4	&	4	&	-4.65 &	 4\\
10	& 124 &	7	&	6	&	-7.41 &	6\\
11	& 174 &	11	&	9	&	-11.34 &5	\\
12	& 167 &	10	&	8	&	-13.27 &4	\\
13	& 138 &	7	&	7	&	-2.77 &	3\\
14	& 55 &	4	&	4	&	-3.54 &	2\\
15	& 312 &	21	&	18	&	-10.38 & 3	\\
Avg. & 151.30 & 8.80&7.50&-8.21&3.70\\ \hline
16	& 232 &	13	&	9	&	-18.77 & 31	\\
17	& 121 &	9	&	7	&	-12.75 & 7	\\
18	& 101 &	8	&	6	&	-13.23 & 8 \\
19	& 112 &	9	&	6	&	-17.21 & 7	\\
20	& 97 &	7	&	6	&	-11.43 & 5	\\
21	& 65 &	5	&	3	&	-14.85 & 4	\\
22	& 71 &	6	&	5	&	-7.64 & 2	\\
23	& 275 &	15	&	11	&	-19.66 & 35	\\
24	& 501 &	26	&	20	&	-29.23 & 78	\\
25	& 321 &	16	&	12	&	-14.09 & 44	\\
Avg. & 	189.60 & 11.40 & 8.50 & -15.89 & 22.10\\ \hline

		\end{tabular}
	}
\end{table}

For real world instances we compared the difference in the volume generated by the solutions found by the ILS-ATUPC algorithm and by the company operators. The instances are real-world requests received by the company, the size of the boxes is highly heterogeneous between instances but fairly homogeneous in the instance (i.e., each instance has boxes that are similar with only a low number of different box dimensions, but every instance is different from the other). 
Results are reported in Table \ref{t:reali}. Columns \textit{TU op}, \textit{TU ILS-ATUCP} and \textit{Vol $\Delta$} have the same meaning as before, while columns \textit{$CI^{xy}$ ILS-ATUCP}, column \textit{$CI^{xy}$ Op.}, \textit{$CI^{z}$ ILS-ATUCP} and \textit{$CI^z$ Op.} report the $CI^{xy}$ and $CI^z$ value of the solutions found by ILS-ATUCP and the company's operators, respectively.

\begin{table}[H]
	\centering
	\caption{Comparison on real world instances} \label{t:reali}
	\resizebox*{0.65\textwidth}{!}{
		\begin{tabular}[h!]{|l|c|c|c|c|c|c|c|c|c|c|}
			\hline
			I&N boxes & TU op & TU ILS-ATUCP & Vol $\Delta$ (\%) & $CI^{xy}$ ILS-ATUCP & $CI^{xy}$ Op. & $CI^z$ ILS-ATUCP & $CI^z$ Op. & Time ILS-ATUCP (sec.) \\ \hline
1	& 43 &	5	&	4	&	-5.15 & 0.0913&0.2532 &0.4563 & 0.6754 & 1	\\
2	& 44 &	5	&	4	&	-5.37& 0.1053 &0.2354 &0.4872 & 0.7233	& 3\\
3	& 32 &	3	&	3	&	0&0.0987 &0.2837 &0.5122 & 0.7132 &	4\\
4	& 43 &	3	&	3	&	-3.53& 0.1074 &0.2963 &0.4623 & 0.6983 & 2	\\
5	& 35 &	4	&	4	&	-2.12& 0.0811 &0.2256 &0.4421 & 0.7426 & 3\\
6	& 46 &	5	&	5	&	-1.52& 0.1062&0.2497 &0.4321 & 0.7101 &	 5\\
7	& 62 &	5	&	5	&	0&0.1021 &0.2614 &0.4823 & 0.6872 &	 3\\
8	& 68 &	6	&	7	&	3.49&0.1145 &0.3184 &0.5673 & 0.7123 &	2 \\
9	& 56 &	6	&	5	&	-5.34&0.0966 &0.2538 &0.4084 & 0.7253 &	 1\\
10	& 87 &	11	&	9	&	-10.72& 0.0975 &0.1977 &0.4294 & 0.6873 &3 	\\
11	& 53 &	7	&	6	&	-7.17& 0.0899 &0.3254 &0.4455 & 0.6824 & 4	\\
12	& 55 &	7	&	6	&	-5.63& 01053 &0.2372 &0.5012 & 0.6974 &	2 \\
13	& 69 &	8	&	6	&	-9.46&0.1032 & 0.3122 &0.4536 & 0.7173 & 3	\\
14	& 72 &	9	&	7	&	-12.59&0.1058 &0.2373 &0.4012 & 0.7257 & 5	\\
15	& 84 &	9	&	7	&	-7.32&0.1084 & 0.2345 &0.4522 & 0.6989 & 4	\\
16	& 32 &	3	&	2	&	-4.92&0.0958 & 0.2381 &0.4232 & 0.6735 & 1	\\
17	& 48 &	4	&	4	&	0&0.0993 &0.3153 & 0.4433 & 0.6964 &3	 \\
18	& 128 &	6	&	6	&	0&0.0896 &0.2951 &0.4721 & 0.7171 & 2	 \\
19	& 44 &	5	&	5	&	-4.39&0.1062 &0.2897 &0.4836 & 0.7292 & 4	\\
20	& 46 &	4	&	4	&	-1.03&0.1048 &0.2877 &0.4576 & 0.6934 & 2	\\
21	& 86 &	8	&	8	&	0&0.1101 &0.2085 &0.4509 & 0.6742 & 3	\\
22	& 128 &	8	&	7	&	-6.93&0.1075 &0.1983 &0.4673 & 0.7284 & 5	\\
23	& 156 &	9	&	6	&	-7.55&0.0951 &0.2244 &0.4832 & 0.7453 & 3	\\
24	& 140 &	15	&	13	&	-10.18&0.0897 &0.2542 & 0.4567& 0.6232 & 1	\\
25	& 276 &	21	&	19	&	-7.23&0.0944 &0.1989 &0.4535 & 0.5423 & 2	\\ \hline
Avg. & 77.32	&	7.04	&	6.20	&	-4.5864&0.1002 &0.2573 &0.4609 & 0.6968 & 2.84	\\ \hline

		\end{tabular}
	}
\end{table}

ILS-ATUCP is capable of better consolidating loose packages generating lower volume in 19 instances out of 25. In five instances the operators and ILS-ATUCP use the same number and type of TUs. In  one instance  the operators find a solution using a lower volume. Moreover, the number of TUs used by ILS-ATUCP is lower or equal the one used in operators' solutions for all instances except one. On the other side, the solutions found by the operators present much worse values of both $CI^{xy}$ and $CI^z$ values. The layouts of the solutions provided by the company's operators lead to centers of gravity which are positioned fairly high when compared to the maximum height of the TU, and far from the center of the base of the TU. The solutions found by ILS-ATUCP, on the other side, have a center of gravity that is much more centered on the base and in a  much lower height. This leads to more balanced TUs, increasing the safety of the cargo.

\section{Conclusions}\label{sec:conc}

This paper introduces the Air Transport Unit Consolidation Problem (ATUCP) where loose packages have to be consolidated in TUs before being transport by air. The goal is to determine the best consolidation setting in oder to minimize the cost generated by the number and the volume of the TUS used and, also, to define the best disposal of loose packages inside each TU so that the TU is stable. The problem is extremely complex and finds relevant practical applications in freight forwarding services. Also, the problem is new to the literature.

We proposed a formal description of the problem and an iterated local search heuristic approach based on a construction algorithm for generating an initial feasible solution, followed by local search using different types of neighbourhoods. In addition, a procedure for generating random instances is proposed that allows a comparison with an optimal solutions even for complicated instances.

Computational results on randomly generated instances and on real data show that, on one side, the ILS-ATUCP algorithm produces high-quality solutions when compared with the optimal ones and, on the other side, it largely outperforms the practice of the company providing the real data.

Future research directions includes relaxing assumptions made in the definition of the problem, like non-stackability of TUs. 
Also, from the methodological point of view, an interesting field of research is related to the development of exact solution approaches for the problem.

\bibliography{Letteratura}
\bibliographystyle{apalike}

\section{APPENDIX}

\subsection{Instance generation procedure}\label{sec:instances}

We now describe the procedure to generate  instances for the ATUCP.

In the first phase, each instance is defined by the total volume and weight to be shipped. 100 combinations of the values of total volume V and total wight W are randomly generated and reported in Table \ref{t:istanzone} in the Appendix.he values of V and W are initially generated at random and then adjusted to avoid trivial solutions (for instance, when $W>>>V$ the solution is trivial since the weight capacity of the TU is satisfied with a very low number of boxes and thus the layout problem is negligible). For some instances the values of V and W are specifically chosen in order to study whether ILS-ATUCP  is capable of identifying the best set of TUs needed to solve the instance. We then solve a mathematical model to identify the best set of TUs that cover the given volume and weight. The solution found by the mathematical model is then used as starting point for the box generation procedure. In this phase, three partitioning schemes are used, described in Section \ref{s:part_alg}. Each scheme takes the TUs associated with the solution found by the mathematical model and partition their weight and volume in order to generate boxes. Therefore, from each of the 100 combinations of values of total weight and volume, three instances are generated, one for each partitioning scheme, for a total of 300 instances.

\subsubsection{Mathematical Model}\label{s_matmod}
The mathematical model used as to generate the instances is the following.
Given $TUS$, the set of available types of TUs, the model makes use of the following parameters and variables:
\begin{itemize}
	\item $V$ is the total volume of the boxes to be packed.
	\item $W$ is the total weight of the boxes to be packed.
	\item $v_{tut}$ is the maximum volume capacity of TU of type $tut$.
	\item $w_{tut}$ is the maximum weight capacity of TU of type $tut$.
	\item $x_{tut}$ is an integer, non-negative variable representing the number of TUs of type $tut$ used .
	
\end{itemize}

The formulation of the mathematical model is the following:
\begin{subequations}\label{model:generation}
	\begin{eqnarray}
	\label{OF}&& \min \sum_{tut \in TUS} v_{tut} x_{tut} + \beta\sum_{tut \in TUS} x_{tut} \\
	\label{C1}&& \sum_{tut \in TUS} v_{tut} x_{tut} \geq V\\
	\label{C2}&& \sum_{tut \in TUS} w_{tut} x_{tut} \geq W\\
	\label{C3}&& x_{tu} \in \mathbb{N}^+ \ \ \ \forall \ tut \in TUS
	\end{eqnarray}
\end{subequations}

The objective of \eqref{OF} is the minimization of the volume and number of the TUs. The objective function is similar to the one of the ATUCP, where the term associated with the center of gravity is missing. We assume that all boxes have the same weight density (instances in which boxes have different weight density are presented in Section \ref{sec:res_real}). Given how we generate boxes according to the three partitioning schemes described later, the best disposal of boxes always generates a center of gravity which is perfectly centered on the base and placed exactly at half the height of the TU. Note that, anyway, the ILS-ATUCP  may not be able to find the best boxes disposition. 

The volume is used as an approximation of costs.  Constraints \eqref{C1} and \eqref{C2} ensures that the total volume and weight of the boxes to be packed are respected. The solution of the model is used as the optimal solution to estimate the effectiveness of the ILS-ATUCP in finding good solutions.

In the generated instances, $TUS$ contained the following types of TU: 120x80x130, 120x80x160, 120x100x130, 120x100x160, 120x120x130, 120x120x160. The solutions found are reported in Table \ref{t:soluzioni_ottime}.

\begin{table}[H]
	\centering
	\caption{Values of total volume (V) and weight (W) for instance generation} \label{t:istanzone}
	\resizebox*{0.35\textwidth}{!}{
		\begin{tabular}[h!]{|l|cc|l|cc|l|cc|}
			\hline
			I & V & W  &I & V & W  &I & V & W    \\ \hline
			1 	&	17 	&	783 	&	36 	&	11 	&	3,266 	&	71 	&	24 	&	11,663 	\\
			2 	&	5 	&	1,579 	&	37 	&	1 	&	1,765 	&	72 	&	22 	&	3,816 	\\
			3 	&	13 	&	4,289 	&	38 	&	26 	&	2,738 	&	73 	&	1 	&	4,500 	\\
			4 	&	30 	&	5,076 	&	39 	&	24 	&	1,525 	&	74 	&	11 	&	8,031 	\\
			5 	&	10 	&	3,463 	&	40 	&	25 	&	182 	&	75 	&	25 	&	14,057 	\\
			6 	&	15 	&	3,242 	&	41 	&	18 	&	3,000 	&	76 	&	28 	&	15,980 	\\
			7 	&	18 	&	1,869 	&	42 	&	7 	&	4,081 	&	77 	&	4 	&	3,476 	\\
			8 	&	10 	&	1,519 	&	43 	&	4 	&	566 	&	78 	&	29 	&	1,025 	\\
			9 	&	30 	&	4,329 	&	44 	&	28 	&	2,693 	&	79 	&	17 	&	2,959 	\\
			10 	&	12 	&	2,433 	&	45 	&	12 	&	3,546 	&	80 	&	26 	&	1,695 	\\
			11 	&	27 	&	113 	&	46 	&	19 	&	3,131 	&	81 	&	26 	&	4,933 	\\
			12 	&	19 	&	2,268 	&	47 	&	29 	&	3,708 	&	82 	&	25 	&	4,002 	\\
			13 	&	20 	&	2,199 	&	48 	&	18 	&	1,230 	&	83 	&	3 	&	1,578 	\\
			14 	&	11 	&	4,103 	&	49 	&	28 	&	279 	&	84 	&	18 	&	4,082 	\\
			15 	&	3 	&	4,713 	&	50 	&	20 	&	865 	&	85 	&	30 	&	14,555 	\\
			16 	&	18 	&	2,346 	&	51 	&	23 	&	1,069 	&	86 	&	29 	&	14,766 	\\
			17 	&	23 	&	3,829 	&	52 	&	13 	&	299 	&	87 	&	4 	&	1,417 	\\
			18 	&	3 	&	2,062 	&	53 	&	12 	&	919 	&	88 	&	18 	&	3,307 	\\
			19 	&	27 	&	4,300 	&	54 	&	30 	&	16,452 	&	89 	&	17 	&	4,401 	\\
			20 	&	16 	&	3,941 	&	55 	&	14 	&	2,617 	&	90 	&	8 	&	1,428 	\\
			21 	&	1 	&	3,087 	&	56 	&	2 	&	2,913 	&	91 	&	20 	&	12,843 	\\
			22 	&	15 	&	4,276 	&	57 	&	29 	&	14,500 	&	92 	&	9 	&	4,156 	\\
			23 	&	10 	&	3,962 	&	58 	&	6 	&	1,088 	&	93 	&	15 	&	11,353 	\\
			24 	&	12 	&	4,268 	&	59 	&	1 	&	949 	&	94 	&	20 	&	4,866 	\\
			25 	&	18 	&	4,348 	&	60 	&	5 	&	253 	&	95 	&	26 	&	12,692 	\\
			26 	&	20 	&	589 	&	61 	&	28 	&	3,647 	&	96 	&	22 	&	11,209 	\\
			27 	&	3 	&	1,480 	&	62 	&	16 	&	11,706 	&	97 	&	26 	&	932 	\\
			28 	&	25 	&	4,872 	&	63 	&	25 	&	15,575 	&	98 	&	14 	&	4,888 	\\
			29 	&	3 	&	2,655 	&	64 	&	10 	&	1,643 	&	99 	&	29 	&	11,311 	\\
			30 	&	25 	&	412 	&	65 	&	13 	&	8,389 	&	100 	&	13 	&	5,656 	\\
			31 	&	22 	&	2,579 	&	66 	&	9 	&	5,016 	&		&		&		\\
			32 	&	8 	&	2,806 	&	67 	&	1 	&	1,500 	&		&		&		\\
			33 	&	4 	&	366 	&	68 	&	26 	&	11,823 	&		&		&		\\
			34 	&	6 	&	4,291 	&	69 	&	2 	&	3,193 	&		&		&		\\
			35 	&	16 	&	3,398 	&	70 	&	15 	&	6,144 	&		&		&		\\ \hline
			
		\end{tabular}
	}
\end{table}

\begin{table}[h!]
	\centering
	\caption{Optimal lower bound solutions} \label{t:soluzioni_ottime}
	\resizebox*{0.40\textwidth}{!}{
		\begin{tabular}[h!]{|l|cccccc|r|}
			\hline
			I&120x80x130&120x80x160&120x100x130&120x100x160&120x120x130&120x120x160 & Total \\ \hline
			1	&	1	&	5	&	1	&	1	&	0	&	2	&	10	\\
			2	&	0	&	1	&	1	&	1	&	0	&	0	&	3	\\
			3	&	0	&	1	&	0	&	5	&	1	&	0	&	7	\\
			4	&	1	&	0	&	0	&	14	&	1	&	0	&	16	\\
			5	&	0	&	1	&	1	&	0	&	0	&	3	&	5	\\
			6	&	0	&	0	&	1	&	7	&	0	&	0	&	8	\\
			7	&	0	&	0	&	0	&	0	&	1	&	7	&	8	\\
			8	&	1	&	1	&	1	&	1	&	2	&	0	&	6	\\
			9	&	1	&	0	&	0	&	14	&	1	&	0	&	16	\\
			10	&	2	&	0	&	0	&	3	&	2	&	0	&	7	\\
			11	&	1	&	1	&	1	&	1	&	0	&	9	&	13	\\
			12	&	6	&	0	&	0	&	0	&	0	&	5	&	11	\\
			13	&	1	&	0	&	0	&	4	&	1	&	4	&	10	\\
			14	&	3	&	0	&	1	&	2	&	1	&	0	&	7	\\
			15	&	5	&	0	&	0	&	0	&	0	&	0	&	5	\\
			16	&	7	&	4	&	2	&	0	&	0	&	0	&	13	\\
			17	&	0	&	0	&	1	&	9	&	1	&	1	&	12	\\
			18	&	0	&	1	&	1	&	0	&	0	&	0	&	2	\\
			19	&	1	&	1	&	3	&	2	&	1	&	6	&	14	\\
			20	&	6	&	0	&	3	&	2	&	0	&	0	&	11	\\
			21	&	2	&	0	&	1	&	0	&	0	&	0	&	3	\\
			22	&	0	&	0	&	1	&	7	&	0	&	0	&	8	\\
			23	&	0	&	1	&	1	&	0	&	0	&	3	&	5	\\
			24	&	2	&	0	&	3	&	2	&	0	&	0	&	7	\\
			25	&	2	&	4	&	6	&	0	&	0	&	0	&	12	\\
			26	&	1	&	0	&	0	&	4	&	1	&	4	&	10	\\
			27	&	0	&	2	&	0	&	0	&	0	&	0	&	2	\\
			28	&	5	&	2	&	0	&	6	&	1	&	1	&	15	\\
			29	&	0	&	0	&	1	&	0	&	1	&	0	&	2	\\
			30	&	4	&	2	&	2	&	0	&	0	&	6	&	14	\\
			31	&	1	&	0	&	1	&	10	&	0	&	0	&	12	\\
			32	&	0	&	0	&	0	&	2	&	1	&	1	&	4	\\
			33	&	2	&	1	&	0	&	0	&	0	&	0	&	3	\\
			34	&	0	&	0	&	0	&	0	&	2	&	1	&	3	\\
			35	&	6	&	0	&	3	&	2	&	0	&	0	&	11	\\
			36	&	3	&	0	&	1	&	2	&	1	&	0	&	7	\\
			37	&	2	&	0	&	0	&	0	&	0	&	0	&	2	\\
			38	&	7	&	0	&	0	&	9	&	0	&	0	&	16	\\
			39	&	10	&	0	&	0	&	0	&	0	&	5	&	15	\\
			40	&	4	&	2	&	2	&	0	&	0	&	6	&	14	\\
			41	&	2	&	4	&	6	&	0	&	0	&	0	&	12	\\
			42	&	1	&	1	&	0	&	1	&	0	&	1	&	4	\\
			43	&	2	&	1	&	0	&	0	&	0	&	0	&	3	\\
			44	&	5	&	0	&	3	&	5	&	4	&	0	&	17	\\
			45	&	2	&	0	&	0	&	3	&	2	&	0	&	7	\\
			46	&	0	&	1	&	10	&	0	&	1	&	0	&	12	\\
			47	&	2	&	2	&	1	&	0	&	8	&	3	&	16	\\
			48	&	2	&	0	&	0	&	2	&	5	&	1	&	10	\\
			49	&	4	&	1	&	1	&	1	&	1	&	7	&	15	\\
			50	&	1	&	0	&	0	&	4	&	1	&	4	&	10	\\
			51	&	0	&	10	&	1	&	1	&	1	&	1	&	14	\\
			52	&	0	&	0	&	0	&	1	&	1	&	4	&	6	\\
			53	&	1	&	1	&	0	&	0	&	0	&	4	&	6	\\
			54	&	21	&	0	&	0	&	1	&	1	&	0	&	23	\\
			55	&	2	&	0	&	0	&	6	&	0	&	0	&	8	\\
			56	&	3	&	0	&	0	&	0	&	0	&	0	&	3	\\
			57	&	2	&	2	&	1	&	0	&	8	&	3	&	16	\\
			58	&	3	&	0	&	0	&	0	&	0	&	1	&	4	\\
			59	&	1	&	0	&	0	&	0	&	0	&	0	&	1	\\
			60	&	0	&	1	&	1	&	1	&	0	&	0	&	3	\\
			61	&	5	&	0	&	3	&	5	&	4	&	0	&	17	\\
			62	&	6	&	1	&	3	&	0	&	0	&	1	&	11	\\
			63	&	17	&	0	&	0	&	1	&	1	&	0	&	19	\\
			64	&	1	&	1	&	1	&	1	&	2	&	0	&	6	\\
			65	&	4	&	0	&	0	&	2	&	1	&	1	&	8	\\
			66	&	2	&	2	&	1	&	0	&	1	&	0	&	6	\\
			67	&	0	&	0	&	0	&	0	&	1	&	0	&	1	\\
			68	&	7	&	0	&	0	&	9	&	0	&	0	&	16	\\
			69	&	2	&	0	&	1	&	0	&	0	&	0	&	3	\\
			70	&	0	&	0	&	1	&	7	&	0	&	0	&	8	\\
			71	&	10	&	0	&	0	&	6	&	0	&	0	&	16	\\
			72	&	1	&	0	&	1	&	10	&	0	&	0	&	12	\\
			73	&	0	&	0	&	0	&	0	&	3	&	0	&	3	\\
			74	&	3	&	0	&	1	&	2	&	1	&	0	&	7	\\
			75	&	17	&	0	&	0	&	1	&	1	&	0	&	19	\\
			76	&	4	&	9	&	1	&	3	&	1	&	0	&	18	\\
			77	&	2	&	0	&	0	&	0	&	1	&	0	&	3	\\
			78	&	1	&	13	&	5	&	0	&	0	&	0	&	19	\\
			79	&	3	&	0	&	1	&	1	&	4	&	1	&	10	\\
			80	&	8	&	1	&	0	&	2	&	2	&	3	&	16	\\
			81	&	7	&	0	&	0	&	9	&	0	&	0	&	16	\\
			82	&	5	&	2	&	0	&	6	&	1	&	1	&	15	\\
			83	&	0	&	2	&	0	&	0	&	0	&	0	&	2	\\
			84	&	2	&	4	&	6	&	0	&	0	&	0	&	12	\\
			85	&	1	&	0	&		&	14	&	1	&	0	&	16	\\
			86	&	2	&	2	&	1	&	0	&	8	&	3	&	16	\\
			87	&	2	&	1	&	0	&	0	&	0	&	0	&	3	\\
			88	&	2	&	4	&	6	&	0	&	0	&	0	&	12	\\
			89	&	3	&	0	&	1	&	1	&	4	&	1	&	10	\\
			90	&	0	&	4	&	0	&	0	&	1	&	0	&	5	\\
			91	&	13	&	0	&	0	&	1	&	1	&	0	&	15	\\
			92	&	1	&	2	&	3	&	0	&	0	&	0	&	6	\\
			93	&	8	&	1	&	1	&	1	&	0	&	0	&	11	\\
			94	&	1	&	1	&	0	&	8	&	1	&	0	&	11	\\
			95	&	7	&	0	&	0	&	9	&	0	&	0	&	16	\\
			96	&	1	&	0	&	1	&	10	&	0	&	0	&	12	\\
			97	&	8	&	1	&	0	&	2	&	2	&	3	&	16	\\
			98	&	2	&	0	&	0	&	6	&	0	&	0	&	8	\\
			99	&	2	&	2	&	1	&	0	&	8	&	3	&	16	\\
			100	&	0	&	1	&	0	&	5	&	2	&	0	&	8	\\ \hline
			
		\end{tabular}
	}
\end{table}

\begin{table}[h!]
	\centering
	\caption{Partitioning algorithm statistics} \label{t:part_data}
	\resizebox*{0.75\textwidth}{!}{
		\begin{tabular}[h!]{|l|ccc|ccc|l|ccc|ccc|}
			\hline
			& \multicolumn{3}{|c|}{N° Boxes} & \multicolumn{3}{c|}{N° unique dimensions}& & \multicolumn{3}{|c|}{N° Boxes} & \multicolumn{3}{c|}{N° unique dimensions}  \\ \hline
			I&A1&A2&A3&A1&A2&A3&I&A1&A2&A3&A1&A2&A3 \\ \hline
1	&	223	&	771	&	699	&	223	&	744	&	4	&	51	&	295	&	1173	&	955	&	293	&	1115	&	5	\\
2	&	65	&	230	&	244	&	62	&	221	&	3	&	52	&	146	&	601	&	447	&	144	&	586	&	3	\\
3	&	155	&	637	&	607	&	155	&	609	&	4	&	53	&	133	&	593	&	408	&	131	&	572	&	3	\\
4	&	348	&	1386	&	1463	&	344	&	1328	&	5	&	54	&	447	&	1591	&	1335	&	441	&	1518	&	5	\\
5	&	101	&	485	&	364	&	101	&	458	&	4	&	55	&	182	&	680	&	688	&	180	&	663	&	3	\\
6	&	184	&	734	&	756	&	184	&	706	&	3	&	56	&	36	&	105	&	84	&	36	&	103	&	1	\\
7	&	174	&	808	&	567	&	173	&	773	&	3	&	57	&	306	&	1394	&	1044	&	306	&	1343	&	6	\\
8	&	109	&	470	&	426	&	108	&	455	&	6	&	58	&	72	&	333	&	240	&	72	&	323	&	3	\\
9	&	281	&	1427	&	1463	&	278	&	1329	&	4	&	59	&	9	&	63	&	56	&	9	&	62	&	2	\\
10	&	158	&	566	&	526	&	156	&	546	&	5	&	60	&	58	&	233	&	244	&	58	&	224	&	3	\\
11	&	311	&	1233	&	948	&	309	&	1174	&	5	&	61	&	313	&	1370	&	1264	&	311	&	1327	&	6	\\
12	&	217	&	863	&	698	&	215	&	822	&	3	&	62	&	215	&	783	&	724	&	213	&	752	&	5	\\
13	&	204	&	907	&	791	&	203	&	848	&	5	&	63	&	398	&	1247	&	1111	&	395	&	1180	&	5	\\
14	&	146	&	528	&	507	&	145	&	495	&	6	&	64	&	114	&	473	&	426	&	114	&	463	&	6	\\
15	&	43	&	148	&	125	&	43	&	141	&	1	&	65	&	136	&	618	&	551	&	135	&	586	&	5	\\
16	&	260	&	917	&	816	&	259	&	858	&	4	&	66	&	109	&	446	&	387	&	108	&	421	&	6	\\
17	&	236	&	1090	&	1083	&	234	&	1038	&	4	&	67	&	11	&	51	&	32	&	10	&	47	&	1	\\
18	&	29	&	140	&	144	&	28	&	137	&	3	&	68	&	300	&	1239	&	1256	&	295	&	1177	&	3	\\
19	&	287	&	1247	&	1059	&	283	&	1186	&	6	&	69	&	23	&	96	&	67	&	22	&	96	&	1	\\
20	&	200	&	769	&	780	&	198	&	718	&	4	&	70	&	206	&	723	&	756	&	203	&	693	&	2	\\
21	&	14	&	60	&	25	&	13	&	57	&	1	&	71	&	373	&	1184	&	1136	&	369	&	1112	&	3	\\
22	&	186	&	695	&	756	&	185	&	673	&	2	&	72	&	242	&	1037	&	1100	&	242	&	992	&	4	\\
23	&	114	&	467	&	364	&	114	&	444	&	4	&	73	&	17	&	44	&	25	&	17	&	42	&	1	\\
24	&	136	&	515	&	556	&	135	&	494	&	4	&	74	&	129	&	528	&	507	&	128	&	507	&	6	\\
25	&	233	&	916	&	872	&	230	&	868	&	4	&	75	&	326	&	1228	&	1111	&	326	&	1158	&	5	\\
26	&	241	&	925	&	791	&	237	&	887	&	5	&	76	&	321	&	1377	&	1235	&	319	&	1305	&	6	\\
27	&	23	&	159	&	125	&	23	&	151	&	1	&	77	&	40	&	201	&	159	&	40	&	198	&	3	\\
28	&	316	&	1216	&	1119	&	315	&	1168	&	5	&	78	&	334	&	1461	&	1308	&	330	&	1362	&	4	\\
29	&	38	&	143	&	129	&	38	&	142	&	3	&	79	&	199	&	828	&	672	&	199	&	787	&	6	\\
30	&	273	&	1190	&	952	&	270	&	1118	&	5	&	80	&	350	&	1232	&	1046	&	345	&	1161	&	5	\\
31	&	258	&	1063	&	1100	&	254	&	1004	&	4	&	81	&	349	&	1290	&	1256	&	344	&	1226	&	3	\\
32	&	88	&	351	&	327	&	88	&	338	&	3	&	82	&	328	&	1189	&	1119	&	323	&	1124	&	5	\\
33	&	67	&	229	&	173	&	67	&	218	&	2	&	83	&	35	&	151	&	125	&	35	&	146	&	1	\\
34	&	60	&	285	&	193	&	60	&	273	&	2	&	84	&	254	&	900	&	872	&	254	&	852	&	4	\\
35	&	205	&	771	&	780	&	204	&	735	&	4	&	85	&	353	&	1461	&	1463	&	351	&	1383	&	5	\\
36	&	156	&	557	&	507	&	156	&	536	&	6	&	86	&	336	&	1414	&	1044	&	333	&	1328	&	6	\\
37	&	11	&	42	&	21	&	11	&	39	&	1	&	87	&	70	&	207	&	173	&	70	&	196	&	2	\\
38	&	326	&	1270	&	1256	&	323	&	1209	&	3	&	88	&	266	&	881	&	872	&	265	&	824	&	4	\\
39	&	278	&	1161	&	920	&	278	&	1093	&	3	&	89	&	182	&	798	&	672	&	182	&	759	&	6	\\
40	&	285	&	1203	&	952	&	284	&	1134	&	5	&	90	&	103	&	395	&	319	&	102	&	380	&	3	\\
41	&	210	&	702	&	872	&	210	&	686	&	4	&	91	&	271	&	980	&	887	&	270	&	929	&	5	\\
42	&	75	&	259	&	288	&	75	&	255	&	4	&	92	&	101	&	431	&	436	&	101	&	414	&	4	\\
43	&	73	&	148	&	176	&	73	&	148	&	2	&	93	&	213	&	749	&	692	&	210	&	717	&	4	\\
44	&	357	&	1201	&	1264	&	352	&	1180	&	6	&	94	&	218	&	986	&	951	&	213	&	939	&	5	\\
45	&	131	&	555	&	526	&	128	&	540	&	5	&	95	&	308	&	1267	&	1256	&	305	&	1197	&	3	\\
46	&	241	&	832	&	967	&	239	&	809	&	5	&	96	&	269	&	1033	&	1100	&	268	&	995	&	4	\\
47	&	333	&	1354	&	1044	&	332	&	1296	&	6	&	97	&	281	&	1235	&	1044	&	279	&	1161	&	5	\\
48	&	203	&	915	&	691	&	201	&	878	&	5	&	98	&	131	&	690	&	688	&	131	&	668	&	3	\\
49	&	345	&	1367	&	1035	&	341	&	1317	&	6	&	99	&	335	&	1346	&	1044	&	332	&	1281	&	6	\\
50	&	182	&	964	&	791	&	180	&	939	&	5	&	100	&	135	&	654	&	603	&	135	&	628	&	3	\\ \hline

		\end{tabular}
	}
\end{table}

\end{document}